\documentclass{svjour3}
\usepackage{amsmath,amssymb,multirow}
\usepackage{cite,amsrefs}
\usepackage{bbm}
\usepackage{fullpage}
\usepackage{exscale,subfigure,booktabs}
\usepackage{color}
\usepackage{hyperref}
\usepackage{graphicx}
\usepackage[ruled,vlined]{algorithm2e}
\hypersetup{
	colorlinks=true,
	linkcolor=blue, 
	citecolor=red, 
	urlcolor=blue  } 
\usepackage[weather,alpine,misc,geometry]{ifsym}
\newtheorem{condition}{Condition}
\everymath{\displaystyle}

\newcommand{\LPK}{$(\text{LP}_{C,\overline C})$}

\newcommand{\al}{\alpha}
\newcommand{\be}{\beta}

\newcommand{\la}{\lambda}
\newcommand{\de}{\delta}
\newcommand{\eps}{\varepsilon}
\newcommand{\bx}{\bar x}

\newcommand {\R} {\mathbb R}
\newcommand {\N} {\mathbb N}

\newcommand {\dom} {{\rm dom}\,}

\newcommand {\conv} {{\rm conv}\,}

\newcommand{\norm}[1]{\left\Vert#1\right\Vert}
\newcommand{\abs}[1]{\left\vert#1\right\vert}

\newcommand{\ang}[1]{\left\langle #1 \right\rangle}

\newcommand{\AND}{\quad\mbox{and}\quad}
\newcounter{mycount}


\newcommand{\x}[1]{}
\usepackage{exscale,subfigure}
\usepackage{tcolorbox} 
\usepackage{tikz}
\usetikzlibrary{patterns}
\usetikzlibrary{decorations.pathmorphing} 
\usetikzlibrary{matrix} 
\usetikzlibrary{arrows} 
\usetikzlibrary{calc} 
\usetikzlibrary{positioning,fit,calc} 
\usepackage{pgfplots}
\pgfplotsset{compat=1.13}
\usetikzlibrary{calc}
\usepgfplotslibrary{polar}
\pgfplotsset{compat=newest}
\usepgfplotslibrary{fillbetween}
\usetikzlibrary{arrows.meta,intersections}
\usetikzlibrary{matrix}

\title{Cutting plane algorithms for nonlinear binary optimization}
\author{Hoa T. Bui, Qun Lin, Ryan Loxton}
\institute{
	Hoa T. Bui (\Letter\,)
	\at
	ARC Training Centre for Transforming Maintenance through Data Science, and 
School of Electrical Engineering, Computing and Mathematical Sciences, Curtin University, Australia\\
	 \email{hoa.bui@curtin.edu.au}  
	\and 
	Qun Lin
	\at
	School of Electrical Engineering, Computing and Mathematical Sciences, Curtin University, Australia\\
	\email{q.lin@curtin.edu.au}
	\and 
	Ryan Loxton
	\at
	ARC Training Centre for Transforming Maintenance through Data Science, and 
	School of Electrical Engineering, Computing and Mathematical Sciences, Curtin University, Australia\\
	\email{r.loxton@curtin.edu.au}  
}

\date{Received: date / Accepted: date}
\usepackage{tikz}
\usetikzlibrary{patterns}
\usetikzlibrary{decorations.pathmorphing} 
\usetikzlibrary{matrix} 
\usetikzlibrary{arrows} 
\usetikzlibrary{calc} 
\usetikzlibrary{positioning,fit,calc} 
\usepackage{pgfplots}
\usepackage{pst-plot}
\begin{document}
	\maketitle
	\abstract{
	Current state-of-the-art methods for solving discrete optimization problems are usually restricted to convex settings. 
	In this paper, we propose a general approach based on cutting planes for solving nonlinear, possibly nonconvex, binary optimization problems.
	 We provide a rigorous convergence analysis that quantifies the number of iterations required under different conditions. This is different to most other work in discrete optimization where only finite convergence is proved. Moreover, using tools from variational analysis, we provide necessary and sufficient dual optimality conditions.}
	
	
	\section{Introduction}
	
	In this paper, we consider the following maximization problem, referred to as problem~\eqref{NP}:
	\begin{align}
		\tag{NP}\label{NP}\max \quad & f(x)\\
		\label{nc0}	\text{s.t.}\quad& x\in K,\\
		\label{nc1}\quad& g_j(x) \le 0,\quad j = 1,\ldots,m,\\
		\notag \quad & x\in \{0,1\}^n,
	\end{align}
	where $K$ is a non-empty bounded polyhedral set in $\R^n$ ($n\ge 1$) and the functions $f, g_j:\R^n\to \R$ ($j=1,\ldots,m$) are possibly nonlinear.
	Finding global solutions of the nonlinear binary problem \eqref{NP}, and nonlinear discrete optimization problems in general, has a long history of over fifty years. 
	Recently, the most common deterministic approaches for solving \eqref{NP} are \emph{branch-and-bound} (see \cite{land2010automatic,gupta1983nonlinear,wolsey1999integer}), and \emph{outer approximation} (see \cite{duran1986outer,leyffer1993deterministic,yuan1988methode}). Branch-and-bound requires solving relaxed problems to obtain valid upper bounds. Two major challenges arise with this approach: the bounds obtained from solving the relaxations are generally not tight for nonlinear problems, and the relaxations may not even be solvable (e.g., if the relaxations are nonconvex, then there may be no efficient method for finding global solutions).
	In 1986, Duran and Grossmann  \cite{duran1986outer} proposed an \emph{outer approximation} method for solving a special class of Mixed Integer Nonlinear Programming problems (MINLPs) in which the objective and constraint functions are linear with respect to the integer variables. This approach solves a sequence of mixed integer linear problems to achieve optimal solutions in a finite number of iterations. The outer approximation scheme was then extended in \cite{leyffer1993deterministic,yuan1988methode} to a more general class of  MINLPs in which the functions are convex with respect to the integer variables.
	It is mentioned in \cite{lubin2018polyhedral,mahajan2021minotaur} that outer approximation is more efficient than branch-and-bound because it avoids solving nonlinear relaxation problems. The main idea of outer approximation is to approximate the nonlinear components by a collection of linear functions or closed half spaces. This idea comes with the price of restricting the setting to the convex case, because any function (or set) that can be precisely outer approximated by a collection of linear functions (or closed half spaces) must be convex. This restriction motivates the need for a new solution approach for solving Problem~\eqref{NP}.
	
	The \emph{cutting plane method} was first proposed in the 1950s (see \cite{gilmore158linear}), for 
	solving integer linear problems and extended further in \cite{gomory1958outline} to mixed integer linear problems. Then, Kelly \cite{10.2307/2099058} in 1960, and Cheney and Goldstein \cite{CheneyCold} in 1959, independently introduced the cutting plane approach for solving convex programs. This method closely resembles outer approximation, because it also involves solving a sequence of linear problems until the optimal solution is found. However, unlike outer approximation, which approximates the entire feasible region,
	the cutting plane method requires only two key conditions (see Figure~\ref{F1}): 
	\begin{enumerate}
		\item[1.] each cutting plane must remove at least one new infeasible/non-optimal solution; and
		\item[2.] at least one optimal solution must remain after the addition of each cutting plane.
	\end{enumerate}
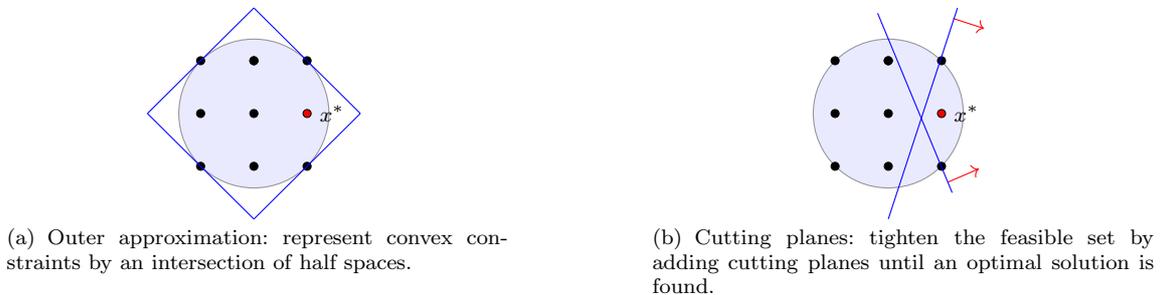
\begin{figure*}[h]
	\centering
	\subfigure[Outer approximation: represent convex constraints by an intersection of half spaces.]{\makebox[6.5cm][c]{
	\begin{tikzpicture}[scale= 0.7]
		\draw[fill = blue!20,opacity=0.4] (0,0) circle (1.41cm);
	\node[circle,draw=black, fill=red, inner sep=0pt,minimum size=3pt, label=right:{$x^*$}] (b) at (1,0) {};
		\node[circle,draw=black, fill = black, inner sep=0pt,minimum size=3pt] (b) at (1,1) {};
			\node[circle,draw=black, fill = black, inner sep=0pt,minimum size=3pt] (b) at (0,1) {};
				\node[circle,draw=black, fill = black, inner sep=0pt,minimum size=3pt] (b) at (1,-1) {};
					\node[circle,draw=black, fill = black, inner sep=0pt,minimum size=3pt] (b) at (0,0) {};
						\node[circle,draw=black, fill = black, inner sep=0pt,minimum size=3pt] (b) at (0,1) {};
							\node[circle,draw=black, fill = black, inner sep=0pt,minimum size=3pt] (b) at (-1,1) {};
								\node[circle,draw=black, fill = black, inner sep=0pt,minimum size=3pt] (b) at (-1,-1) {};
		\node[circle,draw=black, fill = black, inner sep=0pt,minimum size=3pt] (b) at (0,-1) {};
		\node[circle,draw=black, fill = black, inner sep=0pt,minimum size=3pt] (b) at (-1,0) {};
		\draw[blue] (-2,0) -- (0,2);
		\draw[blue] (2,0) -- (0,2);
		\draw[blue] (-2,0) -- (0,-2);
		\draw[blue] (2,0) -- (0,-2);
	\end{tikzpicture}
}
}
\quad \quad \quad \quad\quad\quad
	\subfigure[Cutting planes: tighten the feasible set by adding cutting planes until an optimal solution is found.]{\makebox[6.5cm][c]{
		\begin{tikzpicture}[scale= 0.7]
			\draw[fill = blue!20,opacity=0.4] (0,0) circle (1.41cm);
		\node[circle,draw=black, fill=red, inner sep=0pt,minimum size=3pt, label=right:{$x^*$}] (b) at (1,0) {};
			\node[circle,draw=black, fill = black, inner sep=0pt,minimum size=3pt] (b) at (1,1) {};
			\node[circle,draw=black, fill = black, inner sep=0pt,minimum size=3pt] (b) at (0,1) {};
			\node[circle,draw=black, fill = black, inner sep=0pt,minimum size=3pt] (b) at (1,-1) {};
			\node[circle,draw=black, fill = black, inner sep=0pt,minimum size=3pt] (b) at (0,0) {};
			\node[circle,draw=black, fill = black, inner sep=0pt,minimum size=3pt] (b) at (0,1) {};
			\node[circle,draw=black, fill = black, inner sep=0pt,minimum size=3pt] (b) at (-1,1) {};
			\node[circle,draw=black, fill = black, inner sep=0pt,minimum size=3pt] (b) at (-1,-1) {};
			\node[circle,draw=black, fill = black, inner sep=0pt,minimum size=3pt] (b) at (0,-1) {};
			\node[circle,draw=black, fill = black, inner sep=0pt,minimum size=3pt] (b) at (-1,0) {};
		
			\draw[->,red] (1.12,-1.3)--(1.7,-1.05);
			\draw[->,red] (1.22,1.8)--(1.8,1.61);
			\draw[blue] (0,-2) -- (1.3,2);
			\draw[blue] (-0.2,1.9) -- (1.2,-1.5);
		\end{tikzpicture}
	}
}
	\caption{Applying outer approximation and cutting planes to a circular feasible region.}	\label{F1}
\end{figure*}
	The first condition ensures finite termination, and the second condition ensures convergence to an optimal solution.
	Therefore, the cutting plane method provides flexibility in two respects; {the objective function and constraint set are not required to be convex, and the cutting planes can be, but are not required to be, the tangent planes of the feasible set.} In this paper, we study the cutting plane method for solving the general nonlinear binary problem \eqref{NP}.
	
	We first introduce the framework for our cutting planes approach in Section~\ref{Main}, in which the cutting planes are the tangent planes of the nonlinear functions. 
	Because there is no need for these cutting planes to support the objective function and the feasible set, convexity assumptions can be dropped. However, for convergence to an optimal solution, we require that the region defined by the cutting planes includes at least one optimal solution.
	This is ensured by imposing a special condition (Condition~\ref{con} presented later) that in particular holds when the functions $-f$ and $g(x):=\max_{j=1,\ldots,m} g_j(x)$
	are Lipschitz continuous  and robustly quasiconvex. 
	
	The flexibility of the cutting plane approach becomes apparent when Condition~\ref{con} fails. In such cases, the cutting planes can be modified to ensure that they do not remove all optimal solutions of Problem~\eqref{NP}, and so convergence to an optimal solution still occurs. 
	The idea is to change $v\in \R^n$ and $b\in \R$ in cutting planes of the form $\ang{v,x}\le b$ to ensure Condition~\ref{con} holds.
	In Section~\ref{Application}, we provide two different approaches to achieve this. 
	The first approach introduces the concept of \emph{shifted cutting planes} to modify the component $b$.
	The second approach adds concave and convex penalty terms to the functions $f$ and $g_j$ ($j=1,\ldots,m$) so that Condition~\ref{con} holds; and then the tangent planes of the modified functions can be used as new cutting planes. This approach is essentially modifying the component $v$ (derivatives of $f$ and $g_j$) of the cutting plane $\ang{v,x}\le b$.
	
	For most discrete optimization methods, finite termination is established based on there being only a finite number of points in the search region.
	Therefore, these methods share a common weakness in that, in the worst case, every point is visited (see \cite{hijazi2014outer} for an example on the worst case scenarios of outer approximation). The standard cutting plane algorithm is no exception because if the cutting planes are not tight enough, at each iteration, the generated cuts may only remove one solution from the search region. Accordingly, convergence analysis becomes important to understand computational efficiency and provide guidance on modifying the cutting planes to obtain tighter cuts.
	In Subsection~\ref{estimation}, we show that the tightness of the cutting plane at each iteration is affected by the ratio between the optimality gap (i.e., the difference between the current lower bound and the optimal value of Problem~\eqref{NP}) and the norm of the derivative of the objective function at the current iteration. In other words, modifying the curvature of the objective function can improve the algorithm's computational efficiency. 
	
	The cutting plane method can be viewed as the discrete analogue of \emph{bundle methods}, which have been successfully adapted to solve nonconvex continuous optimization problems. In fact, bundle methods and their variants are currently the most efficient and promising methods for nonsmooth and nonconvex optimization (see \cite{makela2002survey,kiwiel1990proximity}). Like many other algorithms for solving continuous problems, the convergence points of bundle methods are stationary points, i.e., the points at which certain optimality conditions hold. This property is the cornerstone behind the convergence analysis of most algorithms in continuous optimization. 
	It is natural to ponder whether similar characterizations of the cutting plane method exist for binary problems. Such optimality conditions would be useful in two ways. First, when the conditions guaranteeing convergence do not necessarily hold (meaning the algorithm is only heuristic), the optimality conditions can be used  to check whether the final solution is optimal. Second, when the algorithm finds the optimal solution in the first few steps, but the stopping criterion still does not hold (the algorithm may take many more steps to certify optimality), the optimality condition then can be used as an alternative stopping condition.
	 
	The standard optimality conditions are based on certain qualifications that fail in discrete domains. In fact, to the best of our knowledge, there are no tractable optimality conditions for discrete optimization problems.
In Subsection~\ref{Dual}, under the assumption that Condition~\ref{con} holds, we provide dual characterizations for the stopping condition of the cutting plane method. As a by-product of these dual conditions, we show that if the objective function is quasiconvex and the feasible set is defined only with linear constraints, then whenever an optimal solution is encountered, the cutting plane method will converge in the next iteration.
In Subsection~\ref{mod}, we extend the first order Kuhn-Tucker condition from convex programming to the special case of Problem \eqref{NP} in which $-f$ is pseudoconvex and $g_j$ ($j=1,\ldots,m$) are quasiconvex. This condition is useful because it can be used to identify optimal solutions before the cutting plane convergence criteria are satisfied. 
	
	
	The paper is organized as follows. In Section~\ref{Main}, we present the basic methodology of the cutting plane approach for solving nonlinear discrete problems. Under the assumption of Condition~\ref{con}, 
	 Section~\ref{Convergence} contains our convergence results and optimality conditions. 
	 Section~\ref{Application} explains how to apply cutting plane algorithms for solving general nonlinear binary problems. Subsection~\ref{Convex} explores sufficient conditions for when Condition~\ref{con} holds. Furthermore, Subsections~\ref{mod} and \ref{convexification} present two different approaches to deal with the scenarios when 
	 Condition~\ref{con} does not hold. 
	 Finally, in the last section, we test the effectiveness of the approach with an example on quadratic binary problems.
\section{Basic methodology}\label{Main}   

We assume that the functions $f$ and $g_j$ ($j=1,\ldots,m$) are differentiable.
Denote 
\begin{align}\label{setC}
C&:=\left\{x\in K\cap \{0,1\}^n:\; g_j(x)\le 0,\, j =1,\ldots,m \right\},\\
\label{setbC}\overline{C}&:=K\cap \{0,1\}^n\setminus C.
\end{align}
Let $h_f, h_{g_j}:\R^n\times \R^n \to \R$ ($j=1,\ldots,m$) be the
tangent 
planes of $f$ and $g_j$ defined respectively as follows:
\begin{align}
\label{f}	h_f(x,y)&:= \ang{\nabla f(y), x - y}+ f(y),\quad \forall x\in \R^n,\quad \forall y\in \R^n,\\
\label{g1}	h_{g_j}(x,y)&:= \ang{\nabla g_j(y), x - y}+ g_j(y),\quad  \forall x\in \R^n,\quad \forall y\in \R^n, \quad j=1,\ldots,m.
\end{align}
Given two subsets $A_1, A_2\subset \R^n$, let
	\begin{align*}
	\Gamma_{A_1,A_2} :=\bigg\{(x,\theta)\in \R^{n+1}: \;&
	x\in K\cap \{0,1\}^n;\,\theta \le h_f(x,y_1),\; \forall y_1\in A_1;\;\\
	&h_{g_j}(x,y_2)\le 0, \; \forall j \in J(y_2),\;\forall y_2\in A_2
	\bigg\},
	\end{align*}
where $J(y):=\{j:\,g_j(y) = \max_{t=1,\ldots,m}g_t(y)\}$.
Now consider the following auxiliary maximization problem:
\begin{align}
	\label{LP}\tag{$\text{LP}_{A_1,A_2}$}
	\max_{(x,\theta)\in \Gamma_{A_1,A_2}} \;\theta.
\end{align}
This problem can be written explicitly as the following linear binary programming problem:
\begin{align}
\notag	\max \quad &\theta\\
\label{Ax-c1}	\text{s.t.}\quad & \theta \le h_f(x,y), \quad \forall y\in A_1,\\
\label{Ax-c2}	& h_{g_j}(x,y) \le 0,  \quad \forall j\in J(y),\;\forall y\in A_2,\\
\notag	& x\in K\cap \{0,1\}^n.
\end{align}
For $A_1 = C$ and $A_2=\overline C$, the linear problem \LPK\ provides a lower bound for the optimal value of the nonlinear problem \eqref{NP}. This statement is proved in the following proposition.

\begin{proposition} 
	\label{L0} 
The following inequality holds:
	\begin{equation}\label{p1.1}
		\max_{(x,\theta)\in \Gamma_{C,\overline C}} \;\theta \le  \max_{x\in C} f(x).
	\end{equation}
\end{proposition}

\begin{proof}
	Take $(x,\theta)\in \Gamma_{C,\overline C}$.
	We first observe that $x\in C$, since otherwise $x\in \overline C$ and there exists $j\in J(x)$ such that $g_j(x) >0$, and so
	$
	0\ge h_{g_j}(x,x)=g_j(x)>0$,
	which is a contradiction.
	Hence, we must have $x\in C$ and $\theta \le h_f(x,x) = f(x)$.
	This implies $\theta \le f(x) \le \max_{y\in  C} f(y)$, and since $(x,\theta)\in \Gamma_{C,\overline C}$ is arbitrary, inequality \eqref{p1.1} follows.
	\qed
\end{proof}

	The cutting plane algorithm relies on the following condition that requires there is no gap between the optimal values. This is expressed formally in the following condition.
	
	\begin{condition}
		\label{con}
		\begin{equation*}
			\max_{(x,\theta)\in \Gamma_{C,\overline C}} \;\theta = \max_{x\in C} f(x).
		\end{equation*}
	\end{condition}

	We now show that under Condition~\ref{con}, a solution of the nonlinear problem \eqref{NP} can be obtained by solving the linear problem \LPK.	
	\begin{proposition}\label{P2.0}
		Under Condition~\ref{con}, if $(x^*,\theta^*)$ is a solution of \LPK, then $x^*$ is a solution of \eqref{NP}. 
	\end{proposition}
	
	\begin{proof}
	Let $(x^*,\theta^*)$ be a solution of \LPK. 
Then, it follows from the proof of Proposition~\ref{L0} that $x^*$ is feasible for \eqref{NP}, i.e., $x^*\in C$, and by Condition~\ref{con}, 
$$
\max_{x\in C} f(x) = \theta^* \le h_f(x^*,x^*) = f(x^*),
$$
proving the assertion.
		\qed
	\end{proof}
	
	Proposition~\ref{P2.0} also explains that under Condition~\ref{con} the cutting planes \eqref{Ax-c1} and \eqref{Ax-c2} for $A_1= C$ and $A_2 = \overline C$ admit at least one optimal solution of Problem~\eqref{NP}.

We want to find a solution for the nonlinear problem \eqref{NP} by solving the linear problem \LPK. However, it is impossible to generate every cutting plane \eqref{Ax-c1} and \eqref{Ax-c2} in \LPK. Therefore, we propose the following cutting plane algorithm that successively generates cuts of type \eqref{Ax-c1} and \eqref{Ax-c2}.

~

\begin{algorithm}[H]
	\SetAlgoLined
	\vspace{0.5em}
%
	$\eps \ge 0$,
	$\text{UB}_0$ $\leftarrow +\infty $, $\text{LB}_0$ $\leftarrow -\infty$, $k \leftarrow 0$
	
	Take $x^0\in C$
	
	Set $C_1 \leftarrow \{x^0\}$, $\overline C_1 \leftarrow\emptyset$.
	
	\While{  $\text{UB}_{k} - \text{LB}_{k}$ $> \eps$}{
		
		$k \leftarrow k+1$
		
		Solve $(\text{LP}_{C_k, \overline C_k})$ to obtain  $(x^{k},\theta^{k})$

		\eIf{$\exists j\in \{1,\ldots,m\}: g_j(x^k) > 0$}{
			$\overline C_{k+1} \leftarrow \overline C_{k}\cup\{x^k\}$
			
			$\text{UB}_k \leftarrow \theta^k$
		}{
			$C_{k+1} \leftarrow C_{k}\cup\{x^k\}$
			
			$\text{UB}_k \leftarrow \theta^k$, $\text{LB}_k \leftarrow \max\{\text{LB}_{k-1},f(x^k)\}$
			
		}
		
	}
	\label{Al-1}\caption{Cutting plane method for solving \eqref{NP}.}
\end{algorithm}

~

In Algorithm~\ref{Al-1}, the constraints $h_{g_j}(x,x^k) \le 0$ ($j\in J(x^k)$) are added when $x^k$ is not in the feasible region of \eqref{NP} to remove infeasible solutions, in particular $x^k$. These constraints are called \emph{feasibility cuts}. 
Note that we only need to add feasibility cuts corresponding to the constraints $j\in J(x^k)$, i.e., $j\in \{1,\ldots,m\}$ such that $g_j(x^k) =\max_{t=1,\ldots,m}g_t(x^k)$. 

Constraints $\theta \le h_f(x,x^k)$ tighten the optimality gap between the approximation $\min_{y\in C_k} h_f(x^k,y)$ and the value of the objective function $f$, and are called \emph{optimality cuts}. These cuts do not remove the current feasible solution $x^k$, but help to find a better solution.

Consider the sequence  $\{x^k\}$ generated by Algorithm~\ref{Al-1}. 
We have $\Gamma_{C, \overline  C} \subset \Gamma_{C_{k}, \overline C_{k}}$, and the upper bound obtained during iteration $k$ is
$$
\text{UB}_k = \max_{(x,\theta)\in \Gamma_{C_{k},\overline C_{k}}} \theta \ge \max_{(x,\theta)\in \Gamma_{C, \overline C} } \theta.
$$
Moreover, from Condition~\ref{con}, 
$$\text{LB}_k = \max_{x^l\in C_k}f(x^l) \le \max_{x\in C} f(x) = \max_{(x,\theta)\in \Gamma_{C,\overline C}}\theta.$$ Hence,
$$
\text{LB}_k \le \max_{(x,\theta)\in \Gamma_{C,\overline C}} \theta \le \text{UB}_k.
$$
When this expression holds with equality, we have $\max_{(x,\theta)\in \Gamma_{C_k,\overline C_k}}\theta=\max_{(x,\theta)\in \Gamma_{C,\overline C}}\theta=\max_{x\in C}f(x) = f(x^l)$, for some $l\in \{0,1,\ldots,k\}$, and hence $x^l$ is optimal for \eqref{NP}.
We prove in the next theorem that this must occur after a finite number of steps.

%
%

\begin{theorem}\label{T4}
	Under Condition~\ref{con}, the sequence $\{x^k\}\subset K$ generated by Algorithm~\ref{Al-1} (for $\eps =0$) converges to the optimal solution of \eqref{NP} after a finite number of steps.
\end{theorem}

\begin{proof}
	Suppose $x^{k_1} = x^{k_2}$ for some $k_2 > k_1\ge 0$.
Then, $(x^{k_2},\theta^{k_2})\in \arg\max_{(x,\theta)\in \Gamma_{C_{k_2}, \overline C_{k_2}}}\theta$. We consider two cases.
	
	{\bf Case 1.} $x^{k_1} = x^{k_2} \in \overline{C}$. In this case, $x^{k_1}\in \overline C_{k_2}$, and the infeasibility cuts $h_{g_j}(x,x^{k_1}) \le 0$ for $j\in J(x^{k_1})$ are included in $(\text{LP}_{C_{k_2},\overline C_{k_2}})$, therefore $h_{g_j}(x^{k_2},x^{k_1}) \le 0$. We also have
	$$
	h_{g_j}(x^{k_2},x^{k_1}) = g_j(x^{k_1})>0,
	$$
	which is a contradiction.
	
	{\bf Case 2.} $x^{k_1}=x^{k_2} \in {C}$. In this case, $x^{k_1}\in C_{k_2}$, and the optimality cut $\theta \le h_f(x,x^{k_1})$ is included in $(\text{LP}_{C_{k_2},\overline C_{k_2}})$, therefore
	$$
	h_f(x^{k_2},x^{k_1}) = f(x^{k_1})\ge \theta^{k_2}.
	$$
	This implies $\text{UB}_{k_2}$ $= \theta^{k_2} \le f(x^{k_1}) = f(x^{k_2}) \le $ $\text{LB}_{k_2}$, and hence Algorithm~\ref{Al-1} terminates with the optimal solution $x^{k_2}$ at step $k_2$.
	
	The above arguments show that Algorithm~\ref{Al-1} can only revisit a previous point if that point is optimal, and as soon as this occurs the algorithm terminates. Since the set $K\cap\{0,1\}^n$ is finite, we must have finite convergence.
\qed\end{proof}


We now consider Algorithm~\ref{Al-1} when problem \eqref{NP} has nonlinear constraints but a linear objective function. 
In this case, $h_f(x,y) = f(x)$, the optimality cuts become $\theta \le f(x)$, and the algorithm will terminate as soon as Problem $(\text{LP}_{C_{k},\overline C_{k}})$  yields a feasible solution to \eqref{NP}, since when this occurs,
\begin{equation*}
\text{LB}_k \ge f(x^k)= h_{f}(x^k,x^0) \ge \theta^k =  \text{UB}_k.
\end{equation*}
Here, $C_k$ is redundant and Problem $(\text{LP}_{C_{k},\overline C_{k}})$  can be reformulated as
\begin{align}
	\tag{$\text{LP1}_{\overline C_k}$}\label{star}	\max \quad &f(x)\\
	\label{Ax-gc2}	\text{s.t.}\quad & h_{g_j}(x,y) \le 0,  \quad \forall j\in J(y),\;\forall y\in \overline C_k,\\
	\notag	& x\in K\cap \{0,1\}^n.
\end{align}
Problem $(\text{LP}_{C_{k},\overline C_{k}})$  is unbounded when $C_k = \emptyset$ (which is why Algorithm~\ref{Al-1} must commence with a feasible point $x^0\in C$), but Problem~\eqref{star} is always bounded.
Hence, we can solve \eqref{star} successively, without starting from a feasible point, and as soon as the solution $x^k$ is feasible, $(x^k,f(x^k))$ is optimal for $(\text{LP}_{C_{k},\overline C_{k}})$ with $C_k = \{x^k\}$ and $x^k$ is also optimal for \eqref{NP} since by Condition~\ref{con},
$$
f(x^k) \le \max_{x\in C} f(x) = \max_{(x,\theta)\in \Gamma_{C, \overline C} }\theta \le \max_{(x,\theta)\in \Gamma_{C_{k},\overline C_{k}}} \theta = f(x^k).
$$
This discussion leads to the following streamlined version of Algorithm~\ref{Al-1} for the case when $f$ is linear.




~

\begin{algorithm}[H]
	\SetAlgoLined
	
	\vspace{0.5em}
%
	$k \leftarrow 0$
	
	Solve $\max_{x\in K\cap \{0,1\}^n} f(x)$ to obtain $x^0$
	
	Set $\overline C_1  \leftarrow \{x^0\}$
	
	\While{  $\max_{j=1,\ldots,m}g_j(x^k)>0$}{
		$k \leftarrow k+1$
				
		Solve $(\text{LP1}_{\overline C_k})$ to obtain  $x^{k}$
		
		$\overline C_{k+1} \leftarrow \overline C_{k}\cup\{x^k\}$
	
	}
	\label{Al-2}\caption{Cutting plane method for solving \eqref{NP} with linear objective.}
\end{algorithm}
~


It is possible to generalize Algorithm~\ref{Al-1} to nonsmooth cases when the functions $f$ and $g_j$ ($j=1,\ldots,m$) are nonsmooth. In this case, the gradients of $f$ and $g_j$ ($j=1,\ldots,m$) in \eqref{f} and \eqref{g1} can be replaced by the Fr\'echet subdifferentials $\partial \varphi(x)$ and limiting subdifferentials $\overline \partial \varphi(x)$ defined as
	\begin{align}\label{frechet}
	\partial \varphi(x)&:=\left\{v\in \R^n:\; \liminf_{y\to x} \frac{\varphi(y)-\varphi(x) - \ang{v,y-x}}{\norm{y-x}} \ge 0\right\},\\
	\notag
	\overline \partial \varphi(x) &:= \limsup_{y\to x}\partial \varphi(y).
\end{align}
Note that for a convex function, both Fr\'echet subdifferential and limiting subdifferential reduce to the convex subdifferential, i.e.,
$$
\partial \varphi(\bx) = \left\{v\in \R^n:\; \ang{v,x-\bx}\le \varphi(x) - \varphi(\bx)\right\},\quad \forall \bx \in \dom \varphi.
$$
 For a concave function, the subdifferential is 
 $$
 \partial \varphi(\bx) = \left\{v\in \R^n:\; \ang{v,x-\bx}\ge \varphi(x) - \varphi(\bx)\right\},\quad \forall \bx \in \dom \varphi.
 $$
 The results presented in this paper can be generalized to nonsmooth cases by replacing gradients by appropriate subdifferentials.

\section{Convergence results}\label{Convergence}

\subsection{Convergence rate}\label{estimation}
Although Theorem~\ref{T4} shows that Algorithm~\ref{Al-1} terminates after a finite number of steps, in the worst case it might exhaust all points in $K$. 
In this section, we explore the question: when each cutting plane is added, how many non-optimal solutions are eliminated?
Before addressing this question, we first present an important result for our analysis.
\begin{proposition}
	\label{P1}
	Suppose Condition~\ref{con} holds, and let $x^k\in K$ ($k\ge0$) be the iterate generated by Algorithm~\ref{Al-1} during the $k${th} iteration. Then, the following two results hold for the subsequent iterations.
	\begin{enumerate}
			\item  If $x^k\in \overline C$, and $g_j$, $j\in J(x^k)$, have Lipschitz continuous gradients with Lipschitz constants $L(g_j)> 0$, then
		\begin{equation}\label{GC4}
			g_j(x^l)\le\tfrac{1}{2} L(g_j)\norm{x^k-x^l}^2,\quad \forall j\in J(x^k),\quad \forall l> k.
		\end{equation}
		\item If $x^k\in C$, then
		\begin{equation}\label{GC3}
			0\le	\max_{x\in C}f(x) - f(x^k)\le \theta^l- f(x^k)\le \ang{\nabla f(x^{k}), x^l -x^k},\quad \forall l > k.
		\end{equation}
	\end{enumerate}
\end{proposition}

\begin{proof}
	\begin{enumerate}
		\item Suppose $x^k\in \overline{C}$ ($k\ge 0$).
		Recall that for any smooth function $\varphi: \R^n \to \R$ (cf. \cite{bertsekas1998network}) whose gradient is Lipschitz continuous with Lipschitz constant $L(\varphi)$, 
		\begin{equation}\label{amir}
			\varphi(x) \le \varphi(y) +\ang{\nabla\varphi(y),x-y} + \tfrac{1}{2}L(\varphi)\norm{x-y}^2,\quad \forall x,y\in \R^n.
		\end{equation}
		Applying this result to the function $g_j$, where $j\in J(x^k)$, gives
		$$
		g_j(x)\le h_{g_j}(x,x^k)+\tfrac{1}{2} L(g_j) \norm{x-x^k}^2,\quad \forall x\in \R^n.
		$$
		Since $x^k\in \overline C_l$, for $l> k$, we have $h_{g_j}(x^l,x^k) \le 0$, and therefore inequality \eqref{GC4} holds.
		\item 	Suppose $x^k\in C$ ($k\ge 0$). Consider $l>k$, and $(x^l,\theta^l)\in \arg\max_{(x,\theta)\in \Gamma_{C_l,\overline{C}_l}}\theta$.
		From Condition~\ref{con} and $x^k\in C_l$, we have
		$$
		\max_{x\in C} f(x)\le \theta^l\le \ang{\nabla f(x^{k}), x^l -x^k} + f(x^k).
		$$
		Therefore, \eqref{GC3} holds true.
			\qed
	\end{enumerate}
\end{proof}

\begin{remark}
	\label{R6}
	Inequality \eqref{GC4} in Proposition~\ref{P1}(i) shows that at any step $k$, if $x^k\in \overline C$, then in subsequent steps the cutting planes \eqref{Ax-c2} exclude not only $x^k$ but all points $x\in K$ at which the value of the function $g_j$ ($j\in J(x^k)$) is bigger than the quadratic term $\tfrac{1}{2}L(g_j)\norm{x-x^k}^2$. 	
	Thus, in general, the Lipschitz constants $L(g_j)$ ($j=1,\ldots,m$) determine the tightness of the feasibility cuts \eqref{Ax-c2}.
\end{remark}	
	We now use inequality \eqref{GC4} to prove the next result stating that if the Lipschitz constant $L(g_j)> 0$ is sufficiently small for $j\in J(x^k)$, then we have $g_j(x^l) \le 0, \forall l>k$, i.e., the nonlinear constraint $g_j(x)\le 0$ is always satisfied for subsequent iterations. 

	\begin{proposition}\label{PP}
		Suppose Condition~\ref{con} holds. Consider the iterate $x^k$ ($k\ge 0$) generated by Algorithm~\ref{Al-1} during the $k${th} iteration, and assume that $x^k\in \overline C$. For any $j\in J(x^k)$, if the Lipschitz constant for $g_j$ satisfies
		\begin{equation}\label{C3.1}
		0<L(g_j) < \left({2\min_{\substack{x\in K\cap\{0,1\}^n\\ {\small g_j(x) >0}}} g_j(x)}\right)\div \left({\max_{\substack{x\in K\cap\{0,1\}^n\setminus\{x^k\}\\ g_j(x) >0}}\norm{x-x^k}^2}\right),
		\end{equation}
		then 
			$
		 g_j(x^l) \le 0$, for all $l> k$.
	\end{proposition}
\begin{proof}
	Because $x^k\in \overline C$, Proposition~\ref{P1} implies inequality \eqref{GC4}. Take $j\in J(x^k)$, and suppose for some $l>k$, we have $g_j(x^l) >0$. If $x^l=x^k$, then from \eqref{GC4}, we have $0< g_j(x^k) = g_j(x^l)\le \tfrac{1}{2}L(g_j)\norm{x^k-x^l}^2 = 0$, which is a contradiction. Hence, we must have $x^l\neq x^k$.
	
	Because $g_j(x^l) >0$, and $x^l\neq x^k$, from \eqref{C3.1} and \eqref{GC4}, it holds that
$$
 g_j(x^l)\le \tfrac{1}{2}L(g_j)\norm{x^l-x^k}^2 \le 	\tfrac{1}{2}L(g_j)\max_{\substack{x\in K\cap\{0,1\}^n\setminus\{x^k\}\\ g_j(x) >0}}\norm{x-x^k}^2 < \min_{\substack{x\in K\cap\{0,1\}^n\\ {\small g_j(x) >0}}} g_j(x) \le  g_j(x^l),
$$
which is also a contradiction. Hence, we must have $g_j(x^l) \le 0$ for all $l>k$. 
\qed
\end{proof}


\begin{remark}
	When the set $\{{x\in K\cap\{0,1\}^n\setminus\{x^k\}:\, g_j(x) >0}\}$ is empty, we follow the convention $\sup\emptyset = 0$, inequality \eqref{C3.1} becomes
	$$
	0< L(g_j) < +\infty,
	$$
	and $g_j$ is always satisfied in subsequent iterations, irrespective Lipschitz constant.
\end{remark}
	We now examine the second part of Proposition~\ref{P1} in detail. Inequality \eqref{GC3} in Proposition~\ref{P1}(ii) explains how the optimality cuts \eqref{Ax-c1} tighten the feasible set of Problem $(\text{LP}_{C_k,\overline C_k})$. In particular,
	the search region for subsequent iterations is restricted to a closed half space
	\begin{equation}\label{key}
		H_k := \left\{x: c_k \le \ang{\nabla f(x^{k}), x -x^k} \right\},
	\end{equation}
	where $c_k:=  \max_{x\in C}f(x) - f(x^k)$ is the difference between the value of $f$ at $x^k$ and the optimal value of $f$ over $C$. Therefore, the search space
	 after step $k$ becomes
	\begin{equation}\label{key2}
		K_k := K\cap \bigcap_{\substack{t=0,\ldots,k\\x^t\in C}} H_t= K\cap \left(\bigcap_{\substack{t=0,\ldots,k\\x^t\in C}} \left\{x: c_t \le \ang{\nabla f(x^{t}), x -x^t} \right\}\right).
	\end{equation}
	
	Recall that all feasible points of \eqref{NP} are vertices of the $n$-cube $[0,1]^n$. Here, we call such vertices $(0,1)$-vectors.
	We show that the constant 
	\begin{equation}\label{mu}
		\de_k:=\begin{cases}
			\tfrac{\max_{x\in C}f(x)  - f(x^k)}{\norm{\nabla f(x^k)}}, &\quad \text{if } \nabla f(x^k)\neq 0;\\
			0,&\quad \text{otherwise},
		\end{cases}
	\end{equation}
	plays a key role in Algorithm~\ref{Al-1}'s convergence speed.
	We first illustrate this with some special cases.
	\begin{enumerate}
		\item[1.] Suppose $\nabla f(x^k) =0$. Then, Proposition~\ref{P1}(ii) shows that Algorithm~\ref{Al-1} converges in the next iteration.
		\item[2.] Suppose $\de_k > \sqrt{n}$. Then, for any $l >k$, inequality \eqref{GC3} implies
		\begin{equation}\label{R4.0}
		\norm{x^l-x^k}\ge \tfrac{\max_{x\in C}f(x)  - f(x^k)}{\norm{\nabla f(x^k)}} >  \sqrt{n}.
		\end{equation}
		Because the dimension of Problem~\eqref{NP} is $n$, the diameters of $K\cap \{0,1\}^n$ and $C$ are no bigger than the diameter of the cube $[0,1]^n$ in dimension $n$, which is $\sqrt n$. Hence, there is no other point in $K\cap\{0,1\}^n$ that satisfies the above inequality, and Algorithm~\ref{Al-1} must have converged at step $k$. 
		\item[3.] Suppose $\sqrt{n}\ge \de_k > \sqrt{n-1}$, i.e.,
		\begin{equation}\label{R4.1}
		\de_k=\tfrac{\max_{x\in C}f(x)  - f(x^k)}{\norm{\nabla f(x^k)}} >  \sqrt{n-1}.
		\end{equation}
		Then, for any $l >k$, inequality \eqref{GC3} implies that
		$$\norm{x^l-x^k}\ge \tfrac{\max_{x\in C}f(x)  - f(x^k)}{\norm{\nabla f(x^k)}} >  \sqrt{n-1}.$$
		There is only one $(0,1)$-vector $x$ in the $n$-cube $[0,1]^n$ that satisfies $\norm{x-x^k} > \sqrt{n-1}$, so by Condition~\ref{con} and Theorem~\ref{T4}, Algorithm~\ref{Al-1} must converge by step $k+1$. 
		\item[4.] Suppose $\sqrt{n-1}\ge \de_k > \sqrt{n-2}$. At step $k$, the search region is restricted by the half space \eqref{key}, which excludes all $(0,1)$-vectors of the $n$-cube $[0,1]^n$ that lie in any $(n-2)$-cube with $x^k$. Let $\bar x^k$ be the unique $(0,1)$-vector in the $n$-cube $[0,1]^n$ that is of distance $\sqrt{n}$ to $x^k$. Then, in the $n$-cube $[0,1]^n$, all vertices that do not belong to any $(n-2)$-cube with $x^k$ must be adjacent to the vertex $\bar x^k$; and there are exactly $n$ such vertices. Therefore, the search space $K$ is now left with at most $n+1$ points as candidate solutions for the subsequent steps $l> k$.
	\end{enumerate}
The next proposition generalizes this argument to estimate how many binary points Algorithm~\ref{Al-1} eliminates when an optimality cut is added.
\begin{theorem}\label{P3}
	Under Condition~\ref{con}, at step $k$ of Algorithm~\ref{Al-1}, if $x^k\in C$ and there is $N\in \{0,1,\ldots,n\}$ such that
	\begin{align}\label{T5.1}
	\de_k= \tfrac{\max_{x\in C}f(x)  - f(x^k)}{\norm{\nabla f(x^k)}} > \sqrt{N},
	\end{align}
	then the optimality cut $\theta \le h_f(x,x^k)$, removes at least 
	$
	\sum_{q=0}^N\binom{n}{q}
	$
	 binary points from the cube $[0,1]^n$, where $\binom{a}{b} = \frac{b!(a-b)!}{a!}$ for all $a,b\in \N$, and $a\ge b$.
\end{theorem}

\begin{proof}
	Consider $M,m\in \N$ with $n\ge m \ge M$.  Furthermore, consider an arbitrary binary $m$-cube $P$ in $\R^n$ such that $P$ contains $x^k$. The number of $(0,1)$-vectors that are within distance $\sqrt M$ to $x^k$ does not depend on the choice of the cube $P$ as long as $x^k \in P$, and we denote this number by $u(M,m)$.
	From \eqref{GC3}, $u(N,n)$ is a lower bound for the number of binary points that the cutting plane \eqref{Ax-c1} removes from the $n$-cube $[0,1]^n$. We first consider two boundary cases.
	\begin{enumerate}
		\item When $M=m$, all binary points $x$ in the $m$-cube satisfies $\norm{x-x^k} \le \sqrt{M}$, and hence $u(m,m) = 2^m$. 
		\item When $M=0$, the only binary point $x$  in the $m$-cube satisfying $\norm{x-x^k} = 0$ is $x^k$ itself. Hence, $u(0,m) = 1$.
	\end{enumerate}
	Now consider $F$, a facet of $P$, namely a $(m-1)$-cube in $\R^n$, that contains $x^k$. Then, in the facet $F$, there are $u(M,m-1)$ number of $(0,1)$-vectors that are within distance $\sqrt{M}$ to $x^k$. Observe that there is exactly one other facet $F^0$ of $P$ that is opposite to $F$, namely
	$$
	F^0\cap F = \emptyset,\AND \conv(F\cup F^0) = P,
	$$
	and the union of all $(0,1)$-vectors in $F$ and $F^0$ is the set of $(0,1)$-vectors in $P$. 
	An $(0,1)$-vector in $F^0$ is within distance $\sqrt{M}$ to $x^k$ if and only if the vector is within distance $\sqrt{M-1}$ to the projection of $x^k$ onto the facet $F^0$, denoted $\pi(x^k)$. Hence, the number of $(0,1)$-vectors in $F^0$ that are within distance $\sqrt{M}$ to $x^k$ equals the number of $(0,1)$-vertices in $F^0$ that are within distance $\sqrt{M-1}$ to $\pi (x^k)$, and the number of such $(0,1)$-vertices is exactly $u(M-1,m-1)$. Altogether, we have
	\begin{equation}\label{P5.P1}
	u(M,m) = u(M,m-1)+u(M-1,m-1),\quad \forall M\le m-1.
	\end{equation}
	Therefore, from Bernoulli's Triangle, $u(M,m)$ is the sum of the first $M+1$ binomial coefficients of the binomial expansion with power $m$ (see \eqref{EXPo} below and \cite{neiter2016links}):
	\begin{equation}\label{P5.P2}
		u(M,m) = \sum_{q=0}^{M}\binom{m}{q},\quad \forall M\le m.
	\end{equation}
	We now prove \eqref{P5.P2} by induction on the dimension $m$. The base case $m=1$ holds because, by points (i) and (ii) above, $u(0,1) = 1$ and $u(1,1) = 2$. Suppose \eqref{P5.P2} holds for $m = \tau\le n-1$. We prove that \eqref{P5.P2} is also true for $m = \tau +1$. 
	Recall the Pascal's formula, 
	$$
	\binom{a}{b} = \binom{a-1}{b}+\binom{a-1}{b-1},\quad \forall a,b\in \N,\; b\le a-1,
	$$
	and the binomial expansion
	\begin{equation}\label{EXPo}
	(x+y)^a =\sum_{b=0}^{a} \binom{a}{b} x^{a-b}y^b,\quad \forall x,y\in \R,\quad \forall a\in \N.
	\end{equation}
	If $M=\tau+1$, then from point (i) above and the binomial expansion with $x=y=1$, we have $$u(\tau+1,\tau+1) =2^{\tau+1} = \sum_{q=0}^{\tau+1}\binom{\tau+1}{q}.$$
	If $M =0$, then from point (ii) above, we have $$u(M,\tau+1) = \binom{\tau+1}{0}=1.$$
	If $0<M<\tau+1$, then by the induction hypothesis, equality \eqref{P5.P1}, and Pascal's formula we have
	\begin{align*}
		u(M,\tau+1) &= u(M,\tau)+u(M-1,\tau)\\
		&= \sum_{q=0}^{M}\binom{\tau}{q} + \sum_{q=0}^{M-1}\binom{\tau}{q}\\
		& = \binom{\tau}{0} + \sum_{q=1}^{M}\left(\binom{\tau}{q} +\binom{\tau}{q-1}\right)\\
		&= \binom{\tau}{0} + \sum_{q=1}^{M}\binom{\tau+1}{q}=  \sum_{q=0}^{M}\binom{\tau+1}{q}.
	\end{align*}
Thus, by induction, we have proved \eqref{P5.P2} for all $m\le n$, and the result follows immediately since $u(N,n)$ is a lower bound for the number of binary points removed from $[0,1]^n$.
	\qed
\end{proof}

Theorem~\ref{P3} provides a lower bound on the number of binary solutions that each optimality cut removes, when inequality \eqref{T5.1} holds. 
In the next result, we provide an upper bound on the number of optimality cuts required for convergence.

\begin{theorem}\label{T2}
Suppose Condition~\ref{con} holds and let $N\in \{1,\ldots,n\}$. Then, Algorithm~\ref{Al-1} has at most $2^{n-N}$ iterations where \eqref{T5.1} holds.
\end{theorem}

\begin{proof}
	We first prove that for any integer numbers $d>0$ and $m=0,1,\ldots,d$, in a binary $d$-cube, any collection $S$ of vertices of the $d$-cube satisfying the condition that for any $u,v\in S$, $u\neq v$
	\begin{equation}\label{T2.P1}
	\norm{v-u} > \sqrt m,
	\end{equation}
	has the cardinality at most $2^{d-m}$, i.e., $|S| \le 2^{d-m}$.
	We prove the assertion by induction on the dimension $d\in \N$, $d\ge m$. For the base case $d=m$, it is trivial that the set $S$ cannot have more than one vertex, so clearly $\abs{S} \le 1$. Suppose the assertion holds for the dimension $d= t \ge m$, we now prove that it is also true for $d=t+1$. Assume, to the contrary, that there is a collection of vertices $S$ of a binary $(t+1)$-cube, denoted $P$, such that \eqref{T2.P1} holds, and $\abs{S}\ge 2^{t+1-m}+1$.
	In the cube $P$, consider two facets $F$ and $F^\prime$ that are opposite to each other, i.e.,
	$$
	F\cap F^\prime,\AND P = \conv (F\cup F^\prime).
	$$
	Then, both $F$ and $F^\prime$ are binary $t$-cubes. Consider $S_1 := S\cap F$ and $S_2:=S\cap F^\prime$. Then, we have
	$$
	S_1\cap S_2=\emptyset, \AND S = S_1\cup S_2.
	$$
	Therefore,
	\begin{equation}\label{T2.P2}
	\abs{S_1}+\abs{S_2} = \abs{S}\ge 2^{t+1-m}+1.
	\end{equation}
	Observe that \eqref{T2.P1} is satisfied for every pair of vertices in $S_1$ and $S_2$. Thus, by the induction hypothesis, we have 
	$
	\abs{S_1}\le 2^{t-m}$, $
	\abs{S_2}\le 2^{t-m}$, which contradicts \eqref{T2.P2}.
	The assertion above is proved.
	
	Consider the binary cube $[0,1]^n$ in $\R^n$. Consider the collection of feasible solutions $S$ that Algorithm~\ref{Al-1} has generated such that \eqref{T5.1} holds.  Then, from \eqref{GC3}, for every pair $u,v\in S$ with $u\neq v$, inequality \eqref{T2.P1} holds for $m=N$. Thus, $\abs{S} \le 2^{n-N}$.
	\qed
\end{proof}

As a by-product of Theorem~\ref{T2}, if condition \eqref{T5.1} holds for some $N\in \{1,\ldots,n\}$ at every step $k$, then Algorithm~\ref{Al-1} requires no more than $2^{n-N}$ iterations for convergence. 
However, in many cases, when $k$ is large, as Algorithm~\ref{Al-1} gets closer to an optimal solution, the constant $\delta_k$ can become small. The good news is that the dimension of the search region $K_k$ for future iterations also reduces. 
The arguments in the proofs of Theorems~\ref{P3} and \ref{T2} also follows if the dimension $n$ is replaced by the dimension of the current search space, $\dim K_k$. This motivates the next result
that is if 
$\de_k^2+1$ is larger than the dimension of the current research region, the algorithm must conclude by step $(k+1)$th.

\begin{proposition}\label{T8}
	Suppose Condition~\ref{con} holds. Suppose that at step $k\ge 0$, $x^k\in C$ and
	\begin{equation}
	\de_k= \tfrac{\max_{x\in C}f(x)  - f(x^k)}{\norm{\nabla f(x^k)}} > \sqrt{\dim K_k-1},
	\end{equation}
	where the set $K_{k}$ is defined in \eqref{key2}. Then Algorithm~\ref{Al-1} converges in at most  $k+1$ steps.
\end{proposition}

\begin{proof}
	From Proposition~\ref{P1}(ii), after the optimality cut $h_f(x,x^k)\ge \theta$ has been added, the future iterations must satisfy 
	$$
	\norm{x-x^k} \ge \de_k = \tfrac{\max_{x\in C}f(x)  - f(x^k)}{\norm{\nabla f(x^k)}} >  \sqrt{\dim K_k-1}.
	$$
	Therefore, the search region $K_k$ contains at most one binary solution for the next iterations.
	Hence, Algorithm~\ref{Al-1} must converge by step $k+1$.
	\qed 
\end{proof}

Proposition~\ref{T8} establishes a connection between the constant $\delta_k$, the dimension of the search region at the $k$th iteration and the convergence of Algorithm~\ref{Al-1} in the next step. In general, if the dimension of the research space reduces substantially at each step, fast convergence can be guaranteed. 
In the next result, we focus on how each optimality cut $\theta \le h_f(x,x^k)$ reduces the dimension of the search region for subsequent iterations.

\begin{proposition}\label{T6}
	Suppose Condition~\ref{con} holds. Consider step $k$ of Algorithm~\ref{Al-1} with $x^k\in C$ and let
	\begin{align}
		\label{T6.0}
	S_1^+ := \left\{i:\; x^k_i=1,\; \nabla f(x^k)_i >0 \right\},\quad S_0^+ := \left\{i:\; x^k_i=0,\; \nabla f(x^k)_i > 0 \right\},\\
	\label{T6.01}
	S_1^- := \left\{i:\; x^k_i=1,\; \nabla f(x^k)_i \le 0 \right\},\quad S_0^- := \left\{i:\; x^k_i=0,\; \nabla f(x^k)_i \le 0 \right\}.
	\end{align}
	Then the following results hold.
	\begin{enumerate}
		\item[1.] If $i^+ \in S_1^+$ and 
		\begin{equation}\label{T6.1}
		\nabla f(x^k)_{i^+} > \sum_{i\in S_0^+}\nabla f(x^k)_i -\sum_{i\in S_1^-}\nabla f(x^k)_i,
		\end{equation}
		then $x^l_{i^+} = 1$ for all $l > k$.
		\item[2.] If $i^- \in S_0^-$ and 
		\begin{align}\label{T6.2}
		\nabla f(x^k)_{i^-} < -\sum_{i\in S_0^+}\nabla f(x^k)_i +\sum_{i\in S_1^-}\nabla f(x^k)_i,
		\end{align}
		then $x^l_{i^-} = 0$ for all $l > k$.
	\end{enumerate}
\end{proposition}
\begin{proof}
From \eqref{GC3}, 
we have
\begin{align}
\label{T6P1}
\ang{\nabla f(x^k),x^l} \ge \ang{\nabla f(x^k),x^k},\quad \forall l \ge k. 
\end{align}
\begin{enumerate}
	\item[1.] Suppose $i^+\in  S_1^+$ and \eqref{T6.1} holds. Suppose, to the contrary of the result, that $x^l_{i^+} = 0$, for some $l > k$. Then,
	\begin{align*}
	\ang{\nabla f(x^k),x^k} &= \sum_{i\in S_1^+} \nabla f(x^k)_i + \sum_{i\in S_1^-} \nabla f(x^k)_i \\
	& \overset{\eqref{T6.1}}{>} \sum_{i\in S_1^+\setminus\{i^+\}} \nabla f(x^k)_i + \sum_{i\in S_1^-} \nabla f(x^k)_i +\sum_{i\in S_0^+}\nabla f(x^k)_i -\sum_{i\in S_1^-}\nabla f(x^k)_i\\
	&=\sum_{i\in S_1^+\setminus\{i^+\}} \nabla f(x^k)_i  +\sum_{i\in S_0^+}\nabla f(x^k)_i\ge \ang{\nabla f(x^k),x^l},
	\end{align*}
which contradicts \eqref{T6P1}. Thus, we must have $x^l_{i^+} = 1$.
\item[2.] Suppose $i^-\in  S_0^-$ and \eqref{T6.2} holds. Suppose, to the contrary of the result, that $x^l_{i^-} = 1$ for some $l > k$. Then,
\begin{align*}
	\ang{\nabla f(x^k),x^k} &= \sum_{i\in S_1^+} \nabla f(x^k)_i + \sum_{i\in S_1^-} \nabla f(x^k)_i \\
	& \overset{\eqref{T6.2}}{>} \sum_{i\in S_1^+} \nabla f(x^k)_i + \sum_{i\in S_1^-} \nabla f(x^k)_i - \left(-\sum_{i\in S_0^+}\nabla f(x^k)_i +\sum_{i\in S_1^-}\nabla f(x^k)_i\right)+\nabla f(x^k)_{i^-} \\
	&=\sum_{i\in S_1^+} \nabla f(x^k)_i  +\sum_{i\in S_0^+}\nabla f(x^k)_i+\nabla f(x^k)_{i^-} \ge \ang{\nabla f(x^k),x^l},
\end{align*}
which contradicts \eqref{T6P1}. Hence, we must have $x^l_{i^-} = 0$.
\end{enumerate}
\qed
\end{proof}

\begin{remark}\label{R4}
Consider the linear function $l_k(x):=\ang{\nabla f(x^k),x}$. Inequality \eqref{T6P1} shows that at future iterations of Algorithm~\ref{Al-1}, the value of $l_k$ cannot drop below the level at $x^k$, i.e., $l_k(x^l) \ge l_k(x^k)$ for all $l>k$. Note that any change in the value of a variable $x_i^k$ ($i\in S_0^-\cup S_1^+$), from $0$ to $1$ for $i\in S_0^-$ or from $1$ to $0$ for $i\in S_1^+$, will not increase the value of $l_k$ and it will actually decrease $l_k$ if $\nabla f(x^k)\neq 0$. Proposition~\ref{T6} says that if such a decrease exceeds the maximum possible increase from the terms corresponding to $i\in S_0^+\cup S_1^-$, then $l_k(x^l) \ge l_k(x^k)$ cannot be  maintained. In other words, these changes cannot occur in subsequent iterations, hence the values of these indices must be fixed to either $0$ or $1$, which leads to a reduction in the dimension of the feasible region.
\end{remark}
%
%
%

\subsection{Dual conditions for convergence}\label{Dual}

In this subsection, we establish a connection between Algorithm~\ref{Al-1} and gradient-based methods in nonlinear continuous optimization. 
Recall that gradient-based methods for solving continuous optimization problems of the form $\max_{x\in \Omega} f(x)$ typically converge to critical points, e.g., $d(\nabla f(x^*), N_\Omega(x^*))=0$, where $N_\Omega(x^*)$ is the normal cone of $\Omega$ at the point $x^*$ (see \eqref{normalcone} below). Furthermore, when the optimization problem is convex, such critical points are global minimizers. 
The convergence results for Algorithm~\ref{Al-1} presented in this section 
can be viewed as analogues to the classical results on critical points in continuous optimization.

Recall the definition of the normal cone $N_\Omega(x)$ of the closed set $\Omega$ at the point $x\in \Omega$:
\begin{equation}\label{normalcone}
	N_\Omega(x):=\left\{v\in \R^n:\; \ang{v,y-x}\le 0,\; y\in \Omega \right\}.
\end{equation}
Note that $0\in N_\Omega(x)$ for all $x\in \Omega$, and $N_\Omega(x)$ is a closed convex cone.
\begin{theorem}
	\label{P8} Suppose Condition~\ref{con} holds. If $x^k\in C$ ($k\ge 0$) and
	\begin{equation}\label{P8.1}
		d\left(\nabla f(x^k),N_{K\cap \{0,1\}^n}(x^k) \right) < \frac{M_1 - M_2}{\sqrt n},
	\end{equation}
where $M_1:= \max_{x\in C}f(x)$ is the optimal value of \eqref{NP} and $M_2:= \max_{\substack{x\in C\\ f(x) < M_1}}f(x)$ is the second best value of \eqref{NP},
then $x^k$ is an optimal solution of \eqref{NP}. Furthermore, if 
	\begin{equation}\label{P8.2}
	d\left(\nabla f(x^k),N_{K\cap \{0,1\}^n}(x^k) \right) =0,
\end{equation}
then Algorithm~\ref{Al-1} converges after the conclusion of the kth step.
\end{theorem}

\begin{proof}
	Suppose inequality \eqref{P8.1} holds, but $x^k$ is not an optimal solution for Problem~\eqref{NP}. Then $f(x^k) \le M_2$. Let $x^*$ be the convergence point obtained after the conclusion of Algorithm~\ref{Al-1}. Because Condition~\ref{con} holds, $x^*$ maximizes \eqref{NP}.
	By Proposition~\ref{P1}(ii), we have
	\begin{equation}\label{P8.P1}
		M_1 - M_2  \le f(x^*) - f(x^k)\le \ang{\nabla f(x^k), x^*-x^k}.
	\end{equation}
	From \eqref{P8.1}, and $x^*\in K\cap \{0,1\}^n$, there exists $v\in N_{K\cap \{0,1\}^n}(x^k)$ such that 
	\begin{align}
		\label{P8.P4}
		\norm{v-\nabla f(x^k)} &< \frac{M_1-M_2}{\sqrt n},\\
			\label{P8.P3}
		\ang{v,x^*-x^k} &\le 0.
	\end{align}
	From \eqref{P8.P1} and \eqref{P8.P3}, we have
	\begin{align}
		\label{P8.P5}
		M_1 - M_2\le \ang{\nabla f(x^k) - v, x^*-x^k}.
	\end{align}
	Using the Cauchy-Schwarz inequality, combining \eqref{P8.P4} and \eqref{P8.P5} gives
	$$
	M_1 - M_2 \le \ang{\nabla f(x^k) - v, x^*-x^k} \le \norm{\nabla f(x^k) - v}\cdot\norm{x^*-x^k}< \frac{M_1-M_2}{\sqrt n}\norm{x^*-x^k}.$$	
	Now, since feasible points of Problem~\eqref{NP} are binary vectors in $\R^n$, we have $\norm{x^*-x^k} \le \sqrt n$, so the inequality above becomes
	$
	M_1 - M_2 <\frac{M_1-M_2}{\sqrt n}\sqrt n$,
	which is a contradiction. Hence, $x^k$ must be  an optimal solution of Problem~\eqref{NP} and $\text{LB}_k = \max_{x\in C}f(x)$.
	
	Now suppose \eqref{P8.2} holds, but Algorithm~\ref{Al-1} has not converged at step $k$. Because the normal cone $N_{K\cap\{0,1\}^n}(x^k)$ is closed, inequality \eqref{P8.2} implies $\nabla f(x^k)\in N_{K\cap\{0,1\}^n}(x^k)$. 
	 Since Algorithm~\ref{Al-1} does not converge at step $k$, then
	 the next iterate $x^{k+1}$ satisfies
	\begin{align}
		\label{P8.P7}
		\text{UB}_{k+1} \le f(x^k)+\ang{\nabla f(x^k), x^{k+1}-x^k}.
	\end{align} 
	Because $x^{k+1}\in K\cap\{0,1\}^n$ and $\nabla f(x^k)\in N_{K\cap\{0,1\}^n}(x^k)$, we have
	\begin{equation*}
		\ang{\nabla f(x^k),x^{k+1}-x^k} \le 0.
	\end{equation*}
	Combining the inequality above with \eqref{P8.P7}, we have
	$$
	\text{UB}_{k+1} \le f(x^k)+\ang{\nabla f(x^k), x^{k+1}-x^k}\le  f(x^k) \le \text{LB}_{k}.
	$$
	Hence, Algorithm~\ref{Al-1} must converge in step $k+1$.
	\qed
\end{proof}

Using the fact that $0\in N_{K\cap\{0,1\}^n}(x^k)$, Theorem~\ref{P8} implies the following corollary.
\begin{corollary}\label{C4}
	Suppose Condition~\ref{con} holds. 
	If $x^k\in C$ and $$\norm{\nabla f(x^k)} < \frac{M_1-M_2}{\sqrt n},$$
	where $M_1$ and $M_2$ are defined in Theorem~\ref{P8},
	 then $x^k$ is an optimal solution for \eqref{NP}. Furthermore, if 
	 $
	 \nabla f(x^k) =0$, then
	  Algorithm~\ref{Al-1} converges to $x^k$ after the conclusion of the kth step.
\end{corollary}

We now show that when the function $f$ is quasiconvex, the dual condition \eqref{P8.2} with $C$ in place of $K\cap \{0,1\}^n$ becomes necessary. Recall that a function $\varphi:\R^n\to \R$ is quasiconvex if $\varphi(\al x+(1-\al)y)\le \max\left\{\varphi(x),\varphi(y) \right\}$ for all $x,y\in \R^n$. 

\begin{proposition}
	\label{P9} Suppose the function $f$ is quasiconvex. Then for any optimal solution $x^*$ of Problem~\eqref{NP}, 
	\begin{equation}\label{P9.2}
		d\left(\nabla f(x^*),N_{C}(x^*) \right) =0.
	\end{equation}
\end{proposition}
\begin{proof}
	For any $x\in \conv C$, by Caratheodory's theorem, there are extreme points $x_0,\ldots,x_n\in C$ of $\conv C$ and $\al_0,\ldots,\al_n\in [0,1]$ such that $\al_0+\al_1+\cdots+\al_n =1$ and
	$
	x=\al_0x_0+\cdots+\al_n x_n
	$.
	Because $f$ is quasiconvex, $f(x) = f(\al_0x_0+\cdots+\al_n x_n)\le \max\left\{f(x_0),\ldots,f(x_n)\right\}$, which proves that the optimal value of $\max_{x\in \conv C} f(x)$ is achieved at some extreme point of $\conv C$, and the extreme point belongs to $C$. Hence, $
	\max_{\conv C}f(x) = \max_{x\in C}f(x)$.
	
	Suppose $x^*$ is optimal for Problem~\eqref{NP}, and hence $x^*$ is also optimal for problem $\max_{x\in \conv C}f(x)$. By the standard optimality condition for continuous domain, we have
	$\nabla f(x^*)\in N_{\conv C}(x^*)$.
	Since, $C\subset \conv C$, we must have 
	$
	N_{\conv C}(x^*) \subset N_C(x^*),
	$
	and therefore, $\nabla f(x^*)\in N_{C}(x^*)$. Hence, \eqref{P9.2} holds.
	\qed
\end{proof}

When Problem~\eqref{NP} only has linear constraints ($C = K\cap \{0,1\}^n$) and the function $f$ is quasiconvex, then optimality conditions \eqref{P8.2} and \eqref{P9.2} coincide. Hence, in this case, condition \eqref{P8.2} is necessary and sufficient for convergence of Algorithm~\ref{Al-1}.
The following corollary states that in this case, once an optimal solution is added to $C_k$, i.e., $x^k$ maximizes \eqref{NP} for $k\ge 0$, Algorithm~\ref{Al-1} must converge in the next step.

\begin{corollary}\label{T9}
	Suppose Condition~\ref{con} holds, the objective function $f$ is quasiconvex, and there are only linear constraints, i.e., $C = K\cap \{0,1\}^n$.
	If $x^k$ is optimal for \eqref{NP}, then Algorithm~\ref{Al-1} converges after the conclusion of the $k$th step.
\end{corollary}

\begin{proof}
Suppose $x^k$ ($k\ge 0$) is optimal for \eqref{NP}. By Proposition~\ref{P9}, $x^k$ satisfies equation \eqref{P9.2}, and thus \eqref{P8.2} holds.
Then, Theorem~\ref{P8} implies that Algorithm~\ref{Al-1} must converge in the next step.
\qed
\end{proof}

One common weakness of most iterative algorithms for mixed integer programming (e.g., outer approximation, Bender's decomposition method) is that even when the algorithm has generated an optimal solution, it may still require a large number of additional iterations to recognize the optimality. With the cutting plane algorithm (Algorithm~\ref{Al-1}), 
Corollary~\ref{T9} proves that in a special setting, it converges the first time an optimal solution is added. 

\section{Applying Algorithm~\ref{Al-1} and extensions}\label{Application}

The justification for Algorithm~\ref{Al-1} relies on Condition~\ref{con}. Thus, two questions immediately arise:
\begin{enumerate}
	\item[1.] when does Condition~\ref{con} hold?
	\item[2.] what to do when Condition~\ref{con} does not hold?
\end{enumerate}
This section examines these questions.
\subsection{Sufficient conditions for Condition~\ref{con}}\label{Cont}
\label{Convex}
Since Condition~\ref{con} cannot be checked directly, we provide some specific cases when it is guaranteed to hold.

\begin{proposition}\label{T1}
	Suppose  \eqref{NP} has an optimal solution $x^*$ such that
	\begin{align}
		\label{C1}
		&h_f(x^*,y_1) \ge f(x^*), \quad \forall y_1\in C,\\
		\label{C11}
		&h_{g_j}(x^*,y_2) \le 0,\quad \forall j\in J(y_2),\quad \forall y_2 \in \overline C.
	\end{align}
	Then, Condition~\ref{con} holds.
\end{proposition}

\begin{proof}
	Let $x^*$ be an optimal solution of \eqref{NP} satisfying \eqref{C1} and \eqref{C11}. 
	Then, by the definition of $\Gamma_{C,\overline C}$, we have $(x^*,f(x^*))\in \Gamma_{C,\overline C}$, and so $(x^*,f(x^*))$ is feasible for \LPK. 
	Hence, by Proposition~\ref{L0},
	$$
	\max_{x\in C}f(x)=f(x^*) \le \max_{(x,\theta)\in \Gamma_{C,\overline{C}}}\theta \le \max_{x\in C}f(x),
	$$
	which establishes Condition~\ref{con}.
	\qed
\end{proof}

Proposition~\ref{T1} is still hard to use directly because it requires knowing an optimal solution $x^*$ for Problem~\eqref{NP}.
Thus, stronger assumptions are required in practice. 
For example, if $f$ is concave and $g_j$ ($j=1,\ldots,m$) are convex, then
\begin{align*}
	f(x)&\le \ang{\nabla f(y), x - y}+ f(y) = h_f(x,y),\quad \forall x\in \R^n,\quad \forall y\in \R^n,\\
	g_j(x)&\ge \ang{\nabla g_j(y), x - y}+ g_j(y)= h_{g_j}(x,y),\quad \forall x\in \R^n,\quad \forall y\in \R^n, \quad  \forall j = 1,\ldots,m,
\end{align*}
and conditions \eqref{C1} and \eqref{C11} hold trivially. 
In fact, these inequalities do not need to hold for all $x,y\in \R^n$, but only for $(x,y)$ in the discrete sets $C\times C$ and $C\times \overline C$, respectively. This leads to the following corollary of Proposition~\ref{T1}.

\begin{corollary}\label{cor}
	Suppose functions $f$ and $g_j$, satisfy
	\begin{align}
		\label{C1-1}	f(x) &\le h_f(x,y),\quad \forall x,y\in C,\\
		\label{C1-2}   0	 &\ge h_{g_j}(x,y),\quad \forall j\in J(y),\quad \forall x\in {C},\forall y\in \overline{C}.
	\end{align}
	Then, Condition~\ref{con} holds.
\end{corollary}

Inequalities \eqref{C1-1} and \eqref{C1-2} are much stronger than \eqref{C1} and \eqref{C11}, which only need to hold at a single optimal point. 
Corollary~\ref{cor} can be weakened further since \eqref{C1-1} only needs to hold for $x,y\in C$ such that $f(x) \ge f(y)$. 

We now show that Condition~\ref{con} also holds under quasiconvexity and Lipschitz assumptions.
Recall that a function $\varphi$ is $\al$-robustly quasiconvex, for some $\al >0$, if $\varphi(x)+\ang{v,x}$ is quasiconvex for all $v\in \R^n$ with $\norm{v}< \al$. 
Note that a function $\varphi:\R^n\to \R$ is convex if and only if $\varphi(x) +\ang{v,x}$ is quasiconvex for all $v\in \R^n$ (cf. \cite{j.-p._1977}). Hence, we can think of quasiconvexity as $\al$-robust quasiconvexity with $\al =0$, and convexity as $\al$-robust quasiconvexity with $\al = +\infty$. 

\begin{proposition}
	\label{t3}
	Let $g(x):=\max_{j=1,\ldots,m} g_j(x)$ for $x\in \R^n$, and suppose that both $f$ and $g$ are Lipschitz continuous. Furthermore, suppose that there exists a Lipschitz constant $\al>0$ for $f$ such that $-f$ is $\al$-robustly quasiconvex, and a Lipschitz constant $\be>0$ for $g$ such that $g$ is $\beta$-robustly quasiconvex.
	Then, Condition~\ref{con} holds.
\end{proposition}

\begin{proof}
	We recall the first order characterization for robustly quasiconvex functions from \cite[Theorem 3.1]{bui2019characterizations} and  \cite[Proposition 3.1]{aussel1994subdifferential}:
	a proper lower semicontinuous, possibly nonsmooth, function $\varphi:\R^n\rightarrow {\mathbb{R}}$ is $\tau$-robustly quasiconvex ($\tau\ge 0$) if and only if for every $x,y\in \R^n$ the following implication holds:
	\begin{align} \label{MT}
		\varphi(x)\leq\varphi(y)\;\Longrightarrow\;\langle v ,x-y\rangle\leq-\min\left\{\tau\|y-x\|,\varphi(y)-\varphi(x)\right\},\quad \forall v\in \partial \varphi(y),
	\end{align}
	where $\partial \varphi(x)$ is the Fr\'echet subdifferential of $\varphi$ at $x$ defined by
	\begin{equation}\label{frechet}
		\partial \varphi(x):=\left\{v\in \R^n:\; \liminf_{y\to x} \frac{\varphi(y)-\varphi(x) - \ang{v,y-x}}{\norm{y-x}} \ge 0\right\}.
	\end{equation}
	If the function $\varphi$ is smooth at $x$, then its Fr\'echet subdifferential at $x$ reduces to its gradient at $x$, i.e., $\partial \varphi(x) =\left\{\nabla \varphi(x) \right\}$. Furthermore, if $\varphi(x):=\max_{j=1,\ldots,m}\psi_j(x)$, then
	$$
	\partial \varphi(x) = \conv\left\{\nabla \psi_j(x):\; \psi_j(x)=\varphi(x),\quad j = 1,\ldots,m\right\}.
	$$
Suppose $x^*$ is an optimal solution of \eqref{NP}. We will show that inequalities \eqref{C1} and \eqref{C11} hold for $x^*$. 
Take $y_1\in C$ and $y_2\in \overline C$. Then,
	\begin{gather}
		\label{t3.1}
		- f(x^*) \le - f(y_1), \AND g(x^*) \le 0 < g(y_2).
	\end{gather}
	Applying the characterization in \eqref{MT} to  $-f$, and exploiting Lipschitz continuity, we have
	\begin{align*}
		\ang{\nabla f(y_1), x^*-y_1} 
		&\ge \min \left\{\al\norm{y_1-x^*}, f(x^*) - f(y_1) \right\}\\
		& = f(x^*) - f(y_1).
	\end{align*}
	Similarly, applying \eqref{MT} to  $g$, and noting that $\partial g(y_2) = \conv\left\{\nabla g_j(y_2),\; \forall j\in J(y_2) \right\}$, we obtain
	\begin{align*}
		\ang{\nabla g_j(y_2), x^*-y_2} &\le -\min \left\{\beta\norm{y_2-x^*}, g(y_2) - g(x^*) \right\}\\
		&=g(x^*) - g(y_2)\le -g(y_2), \quad \forall j\in J(y_2).  
	\end{align*}
	The inequalities above show that $h_f(x^*,y_1) \ge f(x^*)$ and $h_{g_j}(x^*,y_2)\le 0$ for $j\in J(y_2)$, as required under Proposition~\ref{T1}.
	\qed
\end{proof}

\begin{remark}\label{R1}
	The Lipschitz continuity assumptions on $f$ and $g$ in Proposition~\ref{t3} are only needed to derive:
	$$
	\min\left\{\al\norm{y_1-x^*}, f(x^*) - f(y_1) \right\}= f(x^*) - f(y_1),\quad \min \left\{\beta\norm{y_2-x^*}, g(y_2) - g(x^*) \right\} = g(y_2)-g(x^*).
	$$
	Therefore, Lipschitz continuity of $f$ and $g$ can be replaced by the following weaker Lipschitz-type assumptions, which are local conditions corresponding to an optimal solution $x^*$ of \eqref{NP}:
	\begin{align*}
		\inf_{\substack{y\in \R^n,\\f(y)\ge f(x^*) }}\al  \norm{x-y} &\ge  f(x^*) - f(x),\quad \forall x\in K,\\
		\inf_{\substack{y\in \R^n,\\g(y)\le g(x^*) }}\beta \norm{x-y} &\ge g(x) - g(x^*),\quad \forall x\in K.
	\end{align*}
	These assumptions are related to weak-sharp minima, or linear error bound, which are important notions in variational analysis. 
	Furthermore, these local Lipschitz-type assumptions are weaker than Lipschitz continuity. For example, consider the function $g(x) := \max\{\cos(e^x)+x-1,-\abs{x-1}+1\}$. The sublevel set of $g$ at $x^*=0$ is $\{x:\, g(x)\le 0\}= \R_-$, and
	$$
	\inf_{\substack{y\in \R^n,\\g(y)\le g(x^*) }} \abs{x-y} = \max\{x,0\} \ge \max\{\cos(e^x)+x-1,-\abs{x-1}+1\} = g(x),\quad \forall x\in \R.
	$$
	
		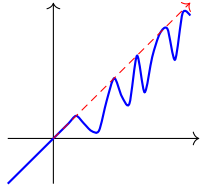
\begin{figure}
		\centering
		\begin{tikzpicture}[scale = 0.6]
			\draw[->] (0, 0) -- (3.2, 0);
			\draw[->] (0, 0) -- (0, 3);
			\draw (-1,0) -- (0,0);
			\draw (0,-1) -- (0,0);
			\draw[thick,scale=0.5, domain=-2:6, smooth, variable=\x, blue] plot ({\x}, {max(\x+cos(deg(2^(\x)))-1, -abs(\x-1)+1)  });
			\draw[->,red, densely dashed] (0, 0) -- (3, 3);
		\end{tikzpicture}
		\caption{Graph of function $g(x) = \max\{\cos(e^x)+x-1,-\abs{x-1}+1\}$.}\label{Fig1}
	\end{figure}

	On the other hand, the function $g$ is not Lipschitz continuous for any constant (see Figure~\ref{Fig1}) since the derivative of $\cos(e^x)+x-1$ is unbounded as $x\to \infty$.
\end{remark}
\begin{remark}
	Even without Lipschitz continuity or the Lipschitz-type conditions in Remark~\ref{R1},  Condition~\ref{con} still holds if $-f$ and $g$ are $\al$- and $\beta$-robustly quasiconvex for sufficiently large $\al$, $\be$, specifically
	$$
	\al \ge \max_{x\in C} f(x) - \min_{x\in C} f(x),\quad \beta \ge \max_{x\in \overline{C}} g(x).
	$$
	In this case, since $x$ and $y$ are binary, $\norm{x-y} \ge 1$ for all $x,y\in K\cap \{0,1\}^n$ with $x\neq y$, and hence,
	\begin{gather*}
		\min\{\al\norm{y-x},f(x) - f(y)\} = f(x) - f(y),\quad \forall x,y\in C,\\
		\min\{\be\norm{y-x},g(y) - g(x)\}\ge 	\min\{\be,g(y)\} = g(y),\quad \forall x\in C,\, y\in \overline{C},
	\end{gather*}
	and the arguments above apply.
\end{remark}

\begin{remark}
	Imposing the conditions in Proposition~\ref{t3} on each individual constraint $g_j(x)\le 0$ ($j=1,\ldots,m$) instead of the maximization $g$ would reduce the power of the result. For example, consider two functions $g_1(x):= x^3$ and $g_2 (x) : = - x^3$. Then, $g(x) = \max\{g_1(x), g_2(x)\} = \abs{x}^3$ is a convex function, and hence $s$-robustly quasiconvex for any $s>0$, but neither $g_1$ or $g_2$ is $s$-robustly quasiconvex for any $s>0$.
\end{remark}

\subsection{Cutting plane modifications for pseudoconvexity}\label{mod}

We now consider the case when Condition~\ref{con} does not hold. 
In this case, when the objective and constraints satisfy pseudoconvexity assumptions, we can modify Algorithm~\ref{Al-1}
by replacing the tangent planes \eqref{f} and \eqref{g1} by the shifted cutting planes, defined as 
\begin{align*}
h_{f}^0(x,y) &:= \ang{\nabla f(y),x-y},\quad \forall x\in \R^n,\quad \forall y\in \R^n\\
h_{g_j}^\eps(x,y) &:= \ang{\nabla g_j(y),x-y}+\eps,\quad \forall x\in \R^n,\quad \forall y\in \R^n,\quad j=1,\ldots,m,
\end{align*}
where $\eps >0$ is a small positive number. We prove later that the modified algorithm, using the shifted cutting planes above instead of \eqref{f} and \eqref{g1}, converges to an optimal solution under pseudoconvexity assumptions. 

Given two subsets $A_1,A_2\subset \R^n$, consider the following auxiliary problem:
\begin{align}
	\tag{$\text{LP2}_{A_1,A_2}$}\label{LP2}	
	\max \quad &\theta\\
	\label{Ax1-c1}	\text{s.t.}\quad & \theta \le h_{f}^0(x,y), \quad \forall y\in A_1,\\
	\label{Ax1-c2}	& h_{g_j}^\eps(x,y) \le 0,  \quad \forall j\in J(y),\;\forall y\in A_2,\\
	\notag	& x\in K\cap \{0,1\}^n.
\end{align}
The optimal value of the linear problem \eqref{LP2} is not necessarily an upper bound for Problem~\eqref{NP}. 
Therefore, the stopping criterion $\text{LB}_k = \text{UB}_k$ in Algorithm~\ref{Al-1} is not applicable with this new auxiliary problem. Note, however, that once Algorithm~\ref{Al-1} repeats a point in $C_k$, no new cutting planes are added and the solution will not change.
Hence, finite termination can be guaranteed if we use the stopping condition $x^k\in C_k$ instead.
We now present the modified algorithm below.

\begin{algorithm}[H]
	\SetAlgoLined
	\vspace{0.5em}
	{\bf Initialization}:
	
	Take $x^0\in C$
	
	$k \leftarrow 1$
	
	Set $C_1 \leftarrow \{x^0\}$, $\overline C_1 \leftarrow \emptyset$

	Solve $(\text{LP2}_{C_1, \overline C_1})$ to obtain  $(x^{1},\theta^{1})$

	\While{  $x^k\notin C_{k}$}{

		\eIf{$\exists j\in \{1,\ldots,m\}: g_j(x^k) > 0$}{
			$\overline C_{k+1} \leftarrow \overline C_{k}\cup\{x^k\}$
			
		}{
			$C_{k+1} \leftarrow C_{k}\cup\{x^k\}$		
		}
	$k \leftarrow k+1$
	
	Solve $(\text{LP2}_{C_k, \overline C_k})$ to obtain  $(x^{k},\theta^{k})$

	}
	\label{Al-3}\caption{modified cutting plane method for solving \eqref{NP}.}
\end{algorithm}

The question is, does the algorithm terminate at an optimal solution?
This depends on whether an optimal solution of \eqref{NP} was added to $C_k$ during a prior iteration.
We prove in the next theorem that
this is indeed the case
 if $-f$ and $g(x):=\max_{j=1,\ldots,m}g_j(x)$ are pseudoconvex. Recall that a function $\varphi:\R^n \to \R$ is pseudoconvex if for any $x,y\in \R^n$,
$$
\varphi(x) < \varphi(y) \Longrightarrow \ang{v,x-y} <0,\quad \text{for all } v\in \partial \varphi(y),
$$
where $\partial \varphi(y)$ is the Fr\'echet subdifferential set of $\varphi$ at $y$.
Pseudoconvexity is stronger than quasiconvexity, but weaker than $\al$-robust quasiconvexity for any $\al >0$. Given a positive number $\al>0$, we have the following implications:
\begin{center}
	Convexity $\Longrightarrow$ $\al$-robust quasiconvexity $\Longrightarrow$ pseudoconvexity $\Longrightarrow$ quasiconvexity.
\end{center}
\begin{theorem}\label{T7}
	Suppose $-f$ and $g(x) :=\max_{t=1,\ldots,m}g_t(x)$ are pseudoconvex. Then, there is a positive number $\bar\eps>0$ such that for any $\eps \in (0,\bar\eps]$, Algorithm~\ref{Al-3} converges to an optimal solution of Problem~\eqref{NP}. 
\end{theorem}

\begin{proof}
	Let $\bar \eps :=  \min_{\substack{(x,y)\in C\times \overline C\\j\in J(y)}} \ang{\nabla g_j(y),y-x}$. For any $x\in C$, $y\in \overline{C}$, by the pseudoconvexity assumption on function $g$, we have
	$$
	g(x)\le 0 < g(y) \Longrightarrow \ang{v,x-y}<0,\quad \text{for all } v\in \partial g(y).
	$$
	Note that $g(y) = \max_{t=1,\ldots,m}g_t(y)$, hence $\nabla g_j(y) \in \partial g(y)$ for all $j\in J(y)$, which implies
	$$
	\ang{\nabla g_j(y),y-x} >0,\quad \forall j\in J(y).
	$$
	Thus, since the sets $C$, $\overline C$ and $J(y)$ are finite, $\bar \eps >0$. For any $\eps \in (0,\bar \eps]$, we have 
	\begin{equation}\label{T3.1}
		h_{g_j}^\eps(x,y) = \ang{\nabla g_j(y),x-y} +\eps \le 0,\quad \forall x\in C,\; y\in \overline{C},\; j\in J(y).
	\end{equation}
	Now, suppose Algorithm~\ref{Al-1} revisits a previous point at step $k\ge 1$, i.e., $x^k\in C_k$, and the set $C_k$ does not contain an optimal solution of \eqref{NP}. Then, constraints \eqref{Ax1-c1} imply that $\theta^k \le \ang{\nabla f(x^k), x^k-x^k} = 0$. Hence, the optimal value of $(\text{LP2}_{C_k, \overline C_k})$ is no larger than $0$.
	
	Let $x^*$ be any solution for \eqref{NP}. Because $C_k$ does not contain any optimal solutions, then $f(y) < f(x^*)$ for all $y\in C_k$.
	From the pseudoconvexity assumption on $-f$, we have $$-f(y) > - f(x^*) \Longrightarrow \ang{-\nabla f(y), x^*-y}<0.$$ 
	Let $\theta :=\min_{y\in C_k} \ang{\nabla f(y), x^*-y}$. Since $C_k$ is finite, we must have $\theta >0$. 
	By the choice of $\theta$, it holds that $(x^*,\theta)$ satisfies constraints \eqref{Ax1-c1}.
	From \eqref{T3.1}, $(x^*,\theta)$ also satisfies constraints \eqref{Ax1-c2}. 
	Hence, $(x^*,\theta)$ is feasible for $(\text{LP2}_{C_k, \overline C_k})$, which implies that the optimal value of $(\text{LP2}_{C_k, \overline C_k})$ is larger than $0$, a contradiction. Thus, $C_k$ must contain an optimal solution of \eqref{NP} after the end of Algorithm~\ref{Al-3}.
	\qed
\end{proof}


To compare the modified cutting planes \eqref{Ax1-c1} with the tangent planes \eqref{Ax-c1}, we first note that the linear problem $(\text{LP}_{A_1,A_2})$ is equivalent to
\begin{align}
	\notag	\max \quad &\theta\\
	\label{Ax-c1'}	\text{s.t.}\quad & \theta  \le h_f(x,y) -M, \quad \forall y\in A_1,\\
	\label{Ax-c2'}	& h_{g_j}(x,y) \le 0,  \quad \forall j\in J(y),\;\forall y\in A_2,\\
	\notag	& x\in K\cap \{0,1\}^n,
\end{align}
where $M:=\max_{x\in C}f(x)$. The cutting planes \eqref{Ax-c1'} are given by
$$
\theta \le \ang{\nabla f(y),x-y}+f(y) - M,\quad \forall x\in \R^n,\quad \forall y\in A_1.
$$
For each $y\in C$, we have $f(y)-M \le 0$, and therefore every pair $(x,\theta)$ that satisfies $\theta  \le h_f(x,y) -M$ also satisfies constraint $\theta \le h_f^0(x,y)$. Hence,
the tangent planes \eqref{Ax-c1'}, which are equivalent to \eqref{Ax-c1}, provide tighter optimality cuts compared to the cutting planes \eqref{Ax1-c1}. 

For the feasibility cuts,
when $\eps$ is sufficiently small (e.g., $\eps <\min_{\substack{j=1,\ldots,m\\ x\in \overline C}}g_j(x)$) the tangent plane \eqref{Ax-c2} implies the cutting plane \eqref{Ax1-c2}, and hence it is also tighter.
The price for the efficiency of the original tangent planes \eqref{Ax-c1}, \eqref{Ax-c2} is that the assumptions required for Condition~\ref{con} to hold, e.g., assumptions on $f$ and $g$ in Proposition~\ref{t3}, are usually stronger than the assumptions in Theorem~\ref{T7}.
%

The pseudoconvexity assumption in Theorem~\ref{T7} explains the motivation for the cuts \eqref{Ax-c1}. When $-f$ is pseudoconvex, if the current point $x^k$ is feasible, but not optimal, then any optimal solution $x^*$ satisfies $f(x^*)> f(x^k)$ and $\ang{\nabla f(x^k),x^*-x^k} >0$. This is consistent with the optimality cut \eqref{Ax1-c1} restricting the search to $\ang{\nabla f(x^k),x-x^k} \ge 0$.

When $g$ is quasiconvex, the constraint $g(u) \le 0$ in \eqref{NP} defines a convex region (even if $g$ itself is not convex), and hence it can be expressed as the intersection of a collection of closed half spaces. 
However, quasiconvexity alone is insufficient, and we require the stronger assumption of pseudoconvexity (which is still much weaker than convexity) to avoid flat level sets, which are often troublesome in quasiconvex optimization (see Example~\ref{EX3} below and  \cite{crouzeix2010geometrical} for further discussions). This is the reason for the pseudoconvexity conditions in Theorem~\ref{T7}.

\begin{figure}
	\centering
	\begin{tikzpicture}[scale = 0.6]
		\draw[->] (0, 0) -- (3.2, 0);
		\draw[->] (0, 0) -- (0, 3);
		\draw (-1,0) -- (0,0);
		\draw (0,-1) -- (0,0);
		\draw[thick,scale=0.5, domain=-1.5:6, smooth, variable=\x, blue] plot ({\x}, {\x+1/4*sin(deg(\x*4)) });
		\node[circle,draw=black, fill=black, inner sep=0pt,minimum size=2pt] (b) at (5/8*pi,0) {};
		\node[circle,draw=black, fill=black, inner sep=0pt,minimum size=2pt] (b) at (1/8*pi,0) {};
		\node[circle,draw=black, fill=black, inner sep=0pt,minimum size=2pt] (b) at (3/8*pi,0) {};
		\node[circle,draw=black, fill=black, inner sep=0pt,minimum size=2pt] (b) at (7/8*pi,0) {};
	\end{tikzpicture}
	\caption{Graph of function $g(x) = \sin(x)+x$.}	\label{Fig2}
\end{figure}
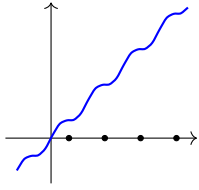

\begin{example}\label{EX3}
	Consider the function $g(x) = \sin(x)+x$ (see Figure~\ref{Fig2}). This function is nondecreasing, and therefore it is quasiconvex. However, $g$ is not pseudoconvex. 
	Consider the set $\overline C=\left\{ (2k+1)\pi: k\ge 0 \right\}$. For any $y \in \overline C$, we have $\nabla g(y) = 0$. Therefore, any cutting plane $h_{g}^\eps$ with $\eps>0$ for $y\in \overline C$ is  invalid:
	$h_{g}^\eps(x,y) = \eps \le 0$.
\end{example}

Because of the finite domain of Problem~\eqref{NP}, even without the conditions in Theorem~\ref{T7}, Algorithm~\ref{Al-3} still terminates, but with no guarantee that the final point is optimal.
 In such cases, we face a new challenge of verifying if the obtained solution is optimal for \eqref{NP}. Without such verification, the algorithm is only heuristic. For convex programming, duality theory provides convenient sufficient conditions to check optimality, but there are no suitable duality conditions 
 for discrete problems. 
Here, we extend the first order Kuhn-Tucker conditions from convex programming to binary problem \eqref{NP} when $-f$ is pseudoconvex, but $g_j$ ($j=1,\ldots,m$) are only required to be quasiconvex. This condition can serve as a sufficient optimality condition to check if the convergence point of Algorithm~\ref{Al-3} is optimal. 
\begin{theorem}\label{P7}
	Suppose $-f$ is pseudoconvex, and $g_j$ ($j=1,\ldots,m$) are quasiconvex. If $x^*\in C$ and there exist constants $\la_j \ge 0 $ such that $\la_jg_j(x^*) = 0$ ($j=1,\ldots,m$) and
	\begin{equation}\label{P7.1}
		\nabla f(x^*)\in N_{K\cap \{0,1\}^n}(x^*) + \la_1 \nabla g_1(x^*)+\cdots+\la_m \nabla g_m(x^*),
	\end{equation}
	then $x^*$ is an optimal solution of Problem~\eqref{NP}. 
\end{theorem}

\begin{proof}
	Suppose there is $\bar x\in C$ such that $f(\bar x) > f(x^*)$. Take $v\in N_{K\cap\{0,1\}^n}(x^*)$ such that 
	\begin{equation}\label{P7.P0}
		\nabla f(x^*) = v+  \la_1\nabla g_1(x^*)+\cdots+\la_m\nabla g_m(x^*).
	\end{equation}
	By the definition of normal cone, 
	\begin{equation}\label{P7.P1}
		\ang{v,\bx - x^*} \le 0.
	\end{equation}
	Because $\la_jg_j(x^*) = 0$ ($j=1,\ldots,m$), then for each $j=1,\ldots,m$, either $\la_j = 0$ or $g_j(x^*) = 0$. Therefore, if $\la_j \neq 0 $, then  $g_j(x^*) =0$, and since $g_j(\bar x) \le 0 = g_j(x^*)$, and $g_j$ is quasiconvex, applying inequality \eqref{MT} with $\tau = 0$ gives
	\begin{align}\label{P7.P2}
		\la_j\ang{\nabla g_j(x^*),\overline x - x^*}\le 0.
	\end{align}
	Note that if $\la_j = 0$, then inequality \eqref{P7.P2} holds trivially. Combining \eqref{P7.P0}, \eqref{P7.P1}, and \eqref{P7.P2}, we have
	\begin{equation}\label{P7.P3}
		\ang{\nabla f(x^*), \bar x - x^*} = \ang{v+  \la_1\nabla g_1(x^*)+\cdots+\la_m\nabla g_m(x^*),\bar x -x^*}\le 0.
	\end{equation}
	On the other hand, $-f$ is pseudoconvex, and hence $-f(\bx) < - f(x^*)$ implies
	$\ang{-\nabla f(x^*), \bar x - x^*} < 0$,
	which contradicts \eqref{P7.P3}. Therefore, $x^*$ must be an optimal solution of $\max_{x\in C} f(x)$.
	\qed
\end{proof}

Computing the normal cone $N_{K\cap \{0,1\}^n}(x^*)$ in \eqref{P7.1} can be difficult
when there is no analytic formulation for the set $K\cap \{0,1\}^n$. 
We show below that condition
\eqref{P7.1} is equivalent to $x^*$ being optimal for the following linear problem 
\begin{equation} \label{C2.1}\tag{LP}
	\max_{x\in K\cap \{0,1\}^n} \ang{\nabla f(x^*) - \la_1 \nabla g_1(x^*)-\cdots-\la_m \nabla g_m(x^*),x},
\end{equation}
which is, compared to \eqref{P7.1}, easier to verify.

\begin{corollary}\label{C2}
	Suppose $-f$ is pseudoconvex, and $g_j$ ($j=1,\ldots,m$) are quasiconvex. If $x^*\in C$ and there exist constants $\la_j\in \R_+$ ($j=1,\ldots,m$) such that $\la_jg_j(x^*) = 0$ ($j=1,\ldots,m$) and $x^*$ is an optimal solution of the linear programming problem \eqref{C2.1},
	then $x^*$ is an optimal solution for \eqref{NP}.
\end{corollary}

\begin{proof}
	Because $x^*\in C$ is the solution for \eqref{C2.1}, then for any $x\in K\cap\{0,1\}^n$, we have 
	\begin{equation}\label{C2.P1}
		\ang{\nabla f(x^*)- \la_1\nabla g_1(x^*)-\cdots-\la_m\nabla g_m(x^*),x - x^*} \le 0,\quad \forall x\in K\cap\{0,1\}^n.
	\end{equation}
	This implies $\nabla f(x^*)- \la_1\nabla g_1(x^*)-\cdots-\la_m\nabla g_m(x^*) \in N_{K\cap \{0,1\}^n}(x^*)$. Therefore, the inclusion \eqref{P7.1} holds. By Theorem~\ref{P7}, $x^*$ is an optimal solution of \eqref{NP}.
	\qed
\end{proof}

\begin{remark}
	When the nonlinear problem \eqref{NP} involves only linear constraints, the multipliers $\la_j$ ($j=1,\ldots,m$) are omitted. Then, Corollary~\ref{C2} implies that if $x^*$ is optimal for $\max_{x\in K\cap \{0,1\}^n}\ang{\nabla f(x^*),x}$, then $x^*$ is also optimal for Problem~\eqref{NP}.
\end{remark}

\subsection{Convexification}\label{convexification}

We now present a second approach to deal with the absence of Condition~\ref{con}. 
Here, instead of modifying the cutting planes as in Subsection~\ref{mod}, we modify the objective and constraints to meet Condition~\ref{con}.
In particular, we exploit the equality constraints 
$x_i^2=x_i$ ($i=1,\ldots,n$), which hold trivially for binary variables, to rewrite the functions $f$ and $g_j$ as $f_\mu(x) := f(x)-\mu\sum_{i=1}^n(x_i^2-x_i)$ and $g_{j,\la_j}(x):=g_j(x)+\la_j\sum_{i=1}^n(x_i^2-x_i)$. When $\mu \ge 0$ and $\la_j\ge 0$ ($j=1,\ldots,m$) are sufficiently large, $f$ and $g_j$ are concave and convex respectively \cite{bertsekas1979convexification}. Following this idea, we focus on finding suitable $\mu$ and $\la_j$ ($j=1,\ldots,m$) so that Condition~\ref{con} is guaranteed.

Consider the following nonlinear binary problem:
\begin{align}
	\tag{$\text{NP}_{\mu,\la_{j}}$}\label{BP}\max \quad & f_{\mu}(x)\\
	\notag\text{s.t.}\quad& x\in K,\\
	\label{nc1-1}\quad& g_{j,\la_{j}}(x) \le 0,\quad j = 1,\ldots,m,\\
	\notag \quad & x\in \{0,1\}^n.
\end{align}
For any choice of $\mu$ and $\la_j$, Problems \eqref{NP} and \eqref{BP} are equivalent.


\begin{lemma}\label{L2}
	Let $\varphi:\R^n\to \R$ be a twice differentiable function and let $\la\ge \tfrac{1}{2}\max\limits_{x\in[0,1]^n\cap K}\lambda_{\max}(\nabla^2\varphi(x))$, where $\lambda_{\max}(\nabla^2\varphi(x))$ is the largest eigenvalue of the Hessian matrix $\nabla^2\varphi(x)$. Then, 
	\begin{equation}\label{la}
		\varphi_\la(x) \le h_{\varphi_\la}(x,y),\quad \forall x,y\in [0,1]^n\cap K,
	\end{equation}
	where $\varphi_\la(x):= \varphi(x)-\la\sum_{i=1}^n(x_i^2-x_i)$ and the cutting plane $h_{\varphi_\la}$ is defined in \eqref{f}.
\end{lemma}

\begin{proof}
	First, we will show that for each $x\in[0,1]^n\cap K$, the Hessian matrix $\nabla^2\varphi_\la(x)$ is negative semidefinite. Indeed,
	$
	\nabla^2\varphi_\la(x)=\nabla^2\varphi(x)-2{\lambda}I_n,
	$
	where $I_n$ is the identity matrix. Let $\xi_1,\ldots,\xi_n$ denote the eigenvalues of $\nabla^2\varphi(x)$. Then for each $i\in\{1,\ldots,n\}$, there exists $z\in \mathbb{R}^n\setminus\{0\}$ such that
	$
	\nabla^2\varphi(x)z= \xi_i z$.
	Hence,
	$$
	\nabla^2\varphi_\la(x)z=\nabla^2\varphi(x)z-2\la z=(\xi_i-2{\lambda})z.
	$$
	The above equality shows that $\xi_1-2{\lambda}, \ldots, \xi_n-2{\lambda}$ are eigenvalues of $\nabla^2\varphi_\la(x)$. Note that 
	$2{\lambda}\ge  \xi_i$, $i=1,\ldots,n$,
	equivalently,
	$
	\xi_i-2{\lambda}\le 0$, $i=1,\ldots,n$.
	Hence, $\nabla^2\varphi_\la(x)$ is a negative semidefinite matrix whenever $x\in[0,1]^n\cap K$.
	
	Let $x$, $y\in [0,1]^n\cap K$. Then, Taylor's theorem implies that 
	$$
	\varphi_\la(x)=\varphi_\la(y)+\langle \nabla \varphi_\la(y),x-y\rangle+\tfrac{1}{2}(x-y)^T\nabla^2\varphi_\la(tx+(1-t)y)(x-y),
	$$
	for some $t\in[0,1]$. Note that since $[0,1]^n\cap K$ is convex, $tx+(1-t)y\in[0,1]^n\cap K$ and hence $\nabla^2\varphi_\la(tx+(1-t)y)$ is negative semidefinite. Thus
	$$
	\varphi_\la(x)-\varphi_\la(y)=\langle \nabla \varphi_\la(y),x-y\rangle+\tfrac{1}{2}(x-y)^T\nabla^2\varphi_\la(tx+(1-t)y)(x-y)\leq \langle \nabla \varphi_\la(y),x-y\rangle.
	$$
	This completes the proof.
	\qed\end{proof}

\begin{proposition}\label{P2}
	Suppose $\mu, \la_{j}\in \R$ ($j=1,\ldots,m$) where 
	$$\mu\ge \tfrac{1}{2}\max\limits_{x\in[0,1]^n\cap  K}\lambda_{\max}(\nabla^2f(x)),\AND \lambda_{j}\ge \tfrac{1}{2}\max\limits_{x\in[0,1]^n\cap K}\lambda_{\max}(\nabla^2g_j(x)).$$
	Then, Condition~\ref{con} holds for Problem~\eqref{BP}.
\end{proposition}

\begin{proof}
	From Lemma~\ref{L2}, and by the choice of $\mu$, $\la_j$ ($j=1,\ldots,m$), we have
	$$
	f_\mu(x)\le h_{f_\mu}(x,y),\AND - g_{j,\la_j}(x) \le -h_{g_{j,\la_j}}(x,y),
	$$
	hold for every $x,y\in K$.
	The first inequality is \eqref{C1-1} with $f_\mu$ in place of $f$.
	The second inequality implies 
	$$
	0\ge g_{j,\la_j}(x) \ge h_{g_{j,\la_j}}(x,y),\quad \forall x\in C,\forall y\in \overline{C},
	$$
	which immediately yields \eqref{C1-2} with $g_{j,\la_j}$ in place of $g_j$. 	
	Hence, the statement in Proposition~\ref{P2} is a direct consequence of Corollary~\ref{cor}.
	\qed
\end{proof}
%
One of the key steps in finding $\mu$ and $\la_{j}$ ($j=1,\ldots,m$) in Proposition~\ref{P2} is to find the maximum eigenvalues for the matrices $\nabla^2 f$ and $\nabla^2 g_j$ ($j=1,\ldots,m$). For matrices with positive entries, the Perron-Frobenius theorem (see \cite{maccluer2000many}) gives the exact value of the largest eigenvalue. A similar result for non-negative matrices is provided below; it shows that the maximum eigenvalue is bounded by the largest row sum.
\begin{lemma}\label{L1}
	Let $A=[a_{ij}]$ be a nonnegative matrix and let $\la_{\max}$ be the largest eigenvalue. Then, 
	$
	\la_{\max} \le \max_{i=1,\ldots,n} \sum_{j=1}^n a_{ij}$. 
\end{lemma}
\begin{proof}
	Observe that $\text{tr}(A) = \sum_{i=1}^n\la_{i}\ge 0$, where $\la_i$ ($i=1,\ldots,n$) are eigenvalues of $A$. Then, $\la_{\max}=\max_{i=1,\ldots,n}\la_i\ge 0$.
	There exists $z\in \R^n\setminus\{0\} $ such that $\max_{i=1,\ldots,n}\abs{z_i} =1$ and $Az=\la_{\max} z$, and hence
	$$
	\la_{\max}\abs{z_i} = \abs{\sum_{j=1}^n a_{ij} z_j},\quad \forall i=1,\ldots,n.
	$$
	Therefore,
	\begin{align*}
		\la_{\max}=\la_{\max}\max_{i=1,\ldots,n}\abs{z_i}  = \max_{i=1,\ldots,n}\abs{\sum_{j=1}^n a_{ij} z_j}&\le \max_{i=1,\ldots,n}\sum_{j=1}^n \abs{a_{ij}z_j}\\
		&= \max_{i=1,\ldots,n}\sum_{j=1}^n \abs{a_{ij}}\abs{z_j}\le \max_{i=1,\ldots,n}\sum_{j=1}^n \abs{a_{ij}}.
	\end{align*}
	\qed\end{proof}

The following example demonstrates the application of Lemma~\ref{L1}, Proposition~\ref{P2} and Algorithm~\ref{Al-1} for solving a binary nonlinear problem.

\begin{example}\label{Ex1}
	Consider the maximization problem
	\begin{align}
		\label{np}\tag{$\text{NP}_1$}	\max\quad &f(x_1,x_2,x_3,x_4)=2x_1x_2x_3+x_1x_3 + 2x_2+3x_3+4x_4\\
		\label{ct1}	\text{s.t.}\quad & 2x_1+x_2+2x_3+2x_4 \le 5\\
		\label{ct2}		& 2x_1+2x_2+x_3+2x_4 \le 5\\
		\label{ct3}		& x_1,x_2,x_3,x_4 \in \{0,1\}.
	\end{align}
We have
	$$
	\nabla f(x_1,x_2,x_3,x_4)= \left( 2x_2x_3+x_3, 2x_1x_3+2,2x_1x_2+x_1+3, 4\right),
	$$
	and 
	$$
	\nabla^2f(x_1,x_2,x_3,x_4) = \begin{bmatrix}
		0& 2x_3 & 2x_2+1&0\\
		2x_3 &  0 & 2x_1 & 0\\
		2x_2+1 & 2x_1 & 0 & 0\\
		0 & 0 & 0 &0
	\end{bmatrix}.
	$$
	By Lemma~\ref{L1}, we set $\mu = 2.5$. Let $f_\mu(x) := f(x)-2.5\sum_{i=1}^3 (x_i^2 - x_i)$. Note that $f$ is linear with respect to $x_4$, and hence the penalty parameter for this variable is irrelevant and has been omitted. The first derivative of the function $f_\mu$ is 
	{
		$$
		\nabla f_\mu(x_1,x_2,x_3,x_4)= \left( 2x_2x_3+x_3-5x_1+2.5, 2x_1x_3-5x_2+4.5,2x_1x_2+x_1-5x_3+5.5, 4\right).
		$$}%
	So now, we can apply Algorithm~\ref{Al-1} to find the maximum of $f_\mu$ over constraints \eqref{ct1}, \eqref{ct2} and \eqref{ct3}. 
	\begin{enumerate}
		\item[1.] Set $x^0 = (1,1,1,0)$, then $\nabla f_\mu(x^0)=(0.5,1.5,3.5,4)$, $f_\mu(x^0)=8$ and solve 
		\begin{align*}
			\label{lp2}\tag{$\text{P}_0$}	\max \quad & \theta\\
			\text{s.t.}\quad & 0.5x_1+1.5x_2+3.5x_3+4x_4 +2.5\ge \theta\\
			& \eqref{ct1},  \eqref{ct2} ,  \eqref{ct3}.\\
		\end{align*}
		We obtain $x^1 =(0,1,1,1)$ as the optimal solution of \eqref{lp2}.
		\item[2.] $\nabla f_\mu(x^1)=(5.5,-0.5,0.5,4)$, $f_\mu(x^1)=9$ and solve 
		\begin{align*}
			\label{lp1}\tag{$\text{P}_1$}	\max \quad & \theta\\
			\text{s.t.}\quad & 0.5x_1+1.5x_2+3.5x_3+4x_4 +2.5\ge \theta\\
			& 5.5x_1-0.5x_2+0.5x_3+4x_4+5\ge \theta\\
			& \eqref{ct1},  \eqref{ct2} ,  \eqref{ct3}.\\
		\end{align*}
		We obtain $x^2 =(0,1,1,1) = x^1$ as the optimal solution of \eqref{lp1}. 
	\end{enumerate}
	Since $(0,1,1,1)$ is a repeated solution, it is the optimal solution of problem \eqref{np}.\qed
\end{example}

\begin{remark}\label{R3}
In Example~\ref{Ex1}, there is no need to add $-3.5(x_4^2 -x_4)$ to $f_\mu$ because the function $f$ is linear with respect to variable $x_4$, and hence $x_4$ has no effect on the Hessian of $f$. The penalty terms $\mu(x_i - x_i^2)$ ($i=1,\ldots,3$) for the other three variables may increase the norm of the gradient $\norm{\nabla f_\mu(x)}$.
Recall from Section~\ref{Convergence} that the ratio $\de_k=\frac{\max_{x\in C}f_\mu(x) - f_\mu(x^k)}{\norm{\nabla f_\mu(x^k)}}$ plays an important role in the speed of Algorithm~\ref{Al-1}. Hence, $\mu$ should be chosen carefully in practice to ensure that both Condition~\ref{con} is satisfied and the effect on $\de_k$ is minimized. Likewise, for problems with nonlinear constraints, the choice of $\la_j$ may effect convergence (see Propositions~\ref{P1}(i) and \ref{PP}).
\end{remark}
%

When the functions $f$ and $g_j$ ($j=1,\ldots,m$) are not twice differentiable or if their second derivatives are computationally expensive, we can use an alternative approach to modify the objective and the constraint functions to meet Condition~\ref{con}. This approach requires that $f$ and $g_j$ ($j=1,\ldots,m$) have Lipschitz continuous gradients.
	Suppose $L(f)$ and $L(g_j)$ ($j=1,\ldots,m$) are the Lipschitz constants of $\nabla f$ and $\nabla g_j$ ($j=1,\ldots,m$), respectively. Then from \eqref{amir}, for all $x\in K\cap\{0,1\}^n$
	\begin{align}
		\label{f_l}	f(x) &= \inf_{y\in K\cap\{0,1\}^n} \left\{h_f(x,y)+\tfrac{1}{2}L(f)\norm{x-y}^2\right\},\\
\label{g_l}	g_j(x) &= \sup_{y\in K\cap\{0,1\}^n}\left\{h_{g_j}(x,y)-\tfrac{1}{2}L(g_j)\norm{x-y}^2\right\},\; j=1,\ldots,m.
	\end{align}
	On the other hand, since $x,y\in \{0,1\}^n$, we have
	\begin{align}
		\notag \norm{x-y}^2 &= \norm{x}^2+\norm{y}^2 -  2\ang{y,x}\\
		\notag 	& = \sum_{i=1}^nx_i-2\ang{y,x}  + \norm{y}^2= \ang{e-2y,x}  + \norm{y}^2,
	\end{align}
	where $e = (1,1,\ldots,1)\in \R^n$. 
	Therefore, for $x\in\{0,1\}^n$, the quadratic functions in the infimum and supremum in \eqref{f_l} and \eqref{g_l}, respectively, can be replaced by the following functions that are linear with respect to $x$:
	\begin{align*}
		h_{f,L(f)}(x,y)&:= h_f(x,y)+\tfrac{1}{2}L(f)\left(\ang{e-2y,x}  + \norm{y}^2\right),\quad \forall x\in \R^n,\quad \forall y\in \R^n,\\
		h_{g_j,L(g_j)}(x,y)&:=h_{g_j}(x,y)-\tfrac{1}{2}L(g_j)\left(\ang{e-2y,x}  + \norm{y}^2\right),\quad  \forall x\in \R^n,\quad \forall y\in \R^n,\; j=1,\ldots,m.
	\end{align*}
From \eqref{f_l} and \eqref{g_l}, we have
\begin{align*}
	f(x) = \inf_{y\in K\cap\{0,1\}^n}h_{f,L(f)}(x,y),\quad g_j(x) = \sup_{y\in K\cap\{0,1\}^n}h_{g_j,L(g_j)}(x,y),\; j=1,\ldots,m,\quad \forall x\in \{0,1\}^n.
\end{align*}
Therefore, Problem~\eqref{NP} is equivalent to the following linear problem:
\begin{align}
\label{NP_L}\tag{$\text{LP}_L$}	\max \quad & \theta\\
\notag\text{s.t.}\quad& x\in K\cap\{0,1\}^n,\\
\notag & \theta \le h_{f,L(f)}(x,y),\quad \forall y\in K\cap \{0,1\}^n\\
\notag\quad& h_{g_j,L(g_j)}(x,y)\le 0, \quad \forall y\in K\cap \{0,1\}^n,\; j\in J(y).
\end{align}
This is an alternative linear problem whose optimal objective value is the same as the original problem, and hence Algorithm~\ref{Al-1} can be applied to the new linear problem. In this formulation, the tangent planes $h_f$ and $h_{g_j}$ are replaced by the planes $h_{f,L(f)}$ and $h_{g_j,L(g_j)}$.
The drawback of this approach is that the Lipschitz constants $L(f)$ and $L(g_j)$ ($j=1,\ldots,m$) might be too large, which can make the cutting planes $h_{f,L(f)}$ and $h_{g_j,L(g_j)}$ inefficient. 

\section{Example: Quadratic Knapsack Problem}
The Quadratic Knapsack Problem (QKP), first introduced in \cite{Gallo1980132}, is a classical binary optimization problem that involves maximizing a quadratic objective subject to a linear capacity constraint (see \cite{PISINGER2007623} for a survey).
We consider the QKP when the weights of all items are identical (see \cite{lima2017solution}). The precise formulation is given below:
\begin{align}
	\tag{QKP0} \label{QKP0}
	\max \quad &f(x):=\tfrac{1}{2}x^TQx+\ang{q,x}\\
	\text{s.t.}\quad & \sum_{i=1}^n x_i \le m,\label{KP}\\
	& x \in \{0,1\}^n,
\end{align}
where $m \ge 2$, $q\in \R_+^n$, and $Q$ is a symmetric $n\times n$ matrix with zero diagonal and positive off-diagonal
entries. 
Since the entries of $q$ and the off-diagonal entries of $Q$ are strictly positive, any optimal solution $x^*$ of \eqref{QKP0} satisfies $\sum_{i=1}^n x_i^* = m \ge 2$. Thus, at least two elements in $x^*$ are non-zero and
\begin{equation}\label{QK-1}
f(x^*) = \tfrac{1}{2}(x^*)^TQ(x^*)+\ang{q,x^*}>\ang{q,x^*}=h_f(x^*,0_{n}),
\end{equation}
where $0_{n}$ is the zero vector in $\R^n$. Therefore, $(x^*,f(x^*))$ does not satisfy the constraint $\theta \le h_f(x,0_{n})$, this implies that Condition~\ref{con} never holds for problem \eqref{QKP0}. Indeed, let $(x^\prime,\theta^\prime)$ be a solution of $(\text{LP}_{C,\overline C})$ so that $\theta^\prime \le h_f(x^\prime,0_n)$.
If Condition~\ref{con} holds, thenby Proposition~\ref{P2.0}, $x^\prime$ is a solution for \eqref{QKP0}, and $\theta^\prime = \max_{x\in C} f(x)=f(x^\prime)$. But this then contradicts \eqref{QK-1}.
Thus Condition~\ref{con} does not hold, and for any solution $x^*$, Algorithm~\ref{Al-1} will never converge to $(x^*,f(x^*))$ when the starting point is the zero vector. Nevertheless,
we now show that, despite Condition~\ref{con} being violated, Algorithm~\ref{Al-1} will converge to $(x^*,f(x^*))$, where $x^*$ is a solution of \eqref{QKP0}, if:
\begin{enumerate}
	\item the initial point $x^0$ belongs in $\tilde C:=\left\{x\in \{0,1\}^n:\; \sum_{i=1}^n x_i =m\right\}$; and 
	\item the matrix $Q$ is conditionally negative definite (c.n.d.), i.e., $x^TQx \le 0$ for any $x\in \R^n$ with $\sum_{i=1}^n x_i = 0$. Note that by \cite[Corollary 4.15]{bapat1997nonnegative}, the matrix $Q$ is c.n.d. if and only if it has exactly one positive eigenvalue. 
\end{enumerate}
Consider a variant of \eqref{QKP0} in which the inequality constraint \eqref{KP} is replaced by the equality constraint $\sum_{i=1}^n x_i = m$:
\begin{align}
	\tag{QKP1} \label{QKP1}
	\max \quad &f(x)\\
	\notag	\text{s.t.}\quad & x\in \tilde{C}.
\end{align}	
This problem is equivalent to Problem \eqref{QKP0} because any optimal  solution of \eqref{QKP0} is optimal for \eqref{QKP1} and vice versa. Furthermore, Condition~\ref{con} holds for Problem \eqref{QKP1} if $Q$ is c.n.d. as we now show. 
For any $x,y\in \tilde C$, we have $\sum_{i=1}^n(x_i-y_i) = \sum_{i=1}^n x_i - \sum_{i=1}^n y_i =0$, and since $Q$ is c.n.d., 
\begin{align*}
	h_f(x,y) - f(x) &= \ang{Qy+q,x-y}+ \tfrac{1}{2}y^TQy+\ang{q,y} - \tfrac{1}{2}x^TQx-\ang{q,x} \\
	& = \ang{Qy,x-y}+ \tfrac{1}{2}(x+y)^TQ(y-x)\\
	& = -\tfrac{1}{2} (x-y)^TQ(x-y) \ge 0,
\end{align*}
which, according to Corollary~\ref{cor}, implies that Condition~\ref{con} holds for \eqref{QKP1}. 
Therefore, by Theorem~\ref{T4}, Algorithm~\ref{Al-1} applied to Problem \eqref{QKP1} converges to an optimal solution.

We now return to Problem \eqref{QKP0} and consider what happens when Algorithm~\ref{Al-1} is initiated with a point in $\tilde{C}$. Since the only constraints in \eqref{QKP0} and \eqref{QKP1} are linear, the linear problem for \eqref{QKP0} at iteration $k$ is ($\text{LP}_{C_k,\emptyset}$) and the corresponding linear problem for \eqref{QKP1}
is  ($\text{LP}_{\tilde C\cap C_k,\emptyset}$). We show that if $x^0\in \tilde C$, then $C_k\subset \tilde{C}$ for all $k\ge 0$.

Let $(x^k,\theta^k)$ be the solution of ($\text{LP}_{C_k,\emptyset}$) and suppose $x^k \notin \tilde C$ for some $k\ge 0$.
Then, from \eqref{KP}, $\sum_{i=1}^nx^k_i < m$. 
Thus, there exists $j\in \{1,\ldots,n\}$ such that $x^k_j =0$. Consider $\hat x\in \{0,1\}^n$ where $\hat x_i = x^k_i$ ($i=1,\ldots,n$, $i\neq j$) and $\hat x_j =1$.
Then, $\sum_{i=1}^n \hat x_i \le m$ and for all $0\le l< k$,
\begin{align*}
	h_f(\hat x,x^l) &=\ang{ Qx^l+q, \hat x - x^l} + f(x^l) \\
	&= \ang{ Qx^l+q, x^k - x^l} + f(x^l) +  (Qx^l+q)_j\\
	&= h_f(x^k,x^l)+(Qx^l+q)_j > h_f(x^k,x^l) \ge \theta^k.
\end{align*}
Choose any $\hat \theta$ satisfying $h_f(\hat x,x^l) > \hat\theta > \theta^k$ for each $l=0,\ldots,k-1$. Then, $(\hat x,\hat \theta)$ is feasible for ($\text{LP}_{C_k,\emptyset}$)  and $\hat \theta > \theta^k$, which contradicts the optimality of $(x^k,\theta^k)$.
Hence, we must have $x^k\in \tilde{C}$ for all $k\ge 0$, and the optimal solution
$(x^k,\theta^k)$ of ($\text{LP}_{C_k,\emptyset}$)  is also the optimal solution of ($\text{LP}_{C_k\cap \tilde{C},\emptyset}$) since the feasible set of ($\text{LP}_{C_k\cap \tilde{C},\emptyset}$) is contained in the feasible set of ($\text{LP}_{C_k,\emptyset}$).

We have proved that applying Algorithm~\ref{Al-1} to \eqref{QKP0}, initiated at a point $x^0\in \tilde{C}$, will generate a sequence of solutions for the linear problem corresponding to \eqref{QKP1}, and this sequence converges to an optimal solution of \eqref{QKP1}, which is also a solution of \eqref{QKP0}. 
This example demonstrates an interesting case where Condition~\ref{con} is not satisfied, but Algorithm~\ref{Al-1}  is still guaranteed to converge to an optimal solution if the starting point is appropriate. 
The reason is that there exists a bounded polyhedron that contains the optimal solution, and when Algorithm~\ref{Al-1} starts inside this polyhedron, all iterations remain in the polyhedron and convergence is guaranteed. 

We now provide numerical results for Algorithm~\ref{Al-1} applied to \eqref{QKP0}.
To generate c.n.d. matrix $Q$, we use the following characterization from \cite[Theorem 4.17]{bapat1997nonnegative}: $Q$ is c.n.d. if and only if $Q$ is the square Euclidean distance matrix of $n$ points,
i.e., $Q= [q_{ij}]$ and $q_{ij} = \norm{v_i - v_j}^2$ for all $i,j = 1,\ldots,n$ and some integer $s\ge 1$ and vectors $v_1,\ldots,v_n\in \R^s$ (see also \cite{Schoenberg35remarksto}).

\begin{figure*}[h]
	\centering
	\subfigure[Average run time]{\makebox[6.5cm][c]{
		\begin{tikzpicture}[scale = 0.5]
			\begin{scope}[local bounding box=plot 1]
				\begin{axis}
					[name=plotter,
					xlabel={n}, ylabel={Average run time (s)}, xmin=0, xmax=1150,
					ymin=0, ymax=480,
					height=10cm,
					width=15cm,
					legend columns=5,
					legend pos=north east,
					axis y line*=left,
					axis x line*=bottom,
					ticklabel style = {font=\large},
					label style = {font=\large},
					xtick={100,200, 300, 400, 500,600,700,800,900,1000,1100},
					ytick={0,50,100,150,200,250,300,350,400,450},
					]
					
					\addplot
					[
					line width=1pt,
					color=blue!60!cyan,
					mark=square,
					]
					coordinates{
						(0,0) (100,4.44) (200,14.2) (300,24.4) (400,42) (500,90.22) (600,125.35) (700,168.56) (800,199) (900,291.17) (1000,305.57) (1100,419.30)};
					
				%
				%
				\end{axis}
			\end{scope}
			\end{tikzpicture}
		}
	}
	\quad \quad \quad \quad\quad\quad
	\subfigure[Average iterations]{\makebox[6.5cm][c]{
			\begin{tikzpicture}[scale = 0.5]
			\begin{scope}[local bounding box=plot 2]
				\begin{axis}[
					name=zoom,
					width=50mm,
					height=50mm,
					at=(plotter.north east),
					axis line style={gray,thick},
					xlabel={{n}}, ylabel={{Average iterations}}, xmin=0, xmax=1110,
					ymin=7.5, ymax=11,
					height=10cm,
					width=15cm,
					legend columns=5,
					legend pos=north east,
					axis y line*=left,
					axis x line*=bottom,
					ticklabel style = {font=\large},
					label style = {font=\large},
					]
					\addplot
					[
					line width=1pt,
					color=cyan!50!gray!70,
					mark=halfsquare*,
					]
					coordinates{
						(100,9.42) (200,9.84) (300,8.16) (400,8.1) (500,9.62) (600,9.8) (700,9.26) (800,8.68) (900,9) (1000,8.3) (1100,9.96)};
				\end{axis}
			\end{scope}
		\end{tikzpicture}
		}
	}
	\caption{Average run time and average iterations of the cutting plane method (for medium-size problems)}\label{fig-100-1000}
\end{figure*}
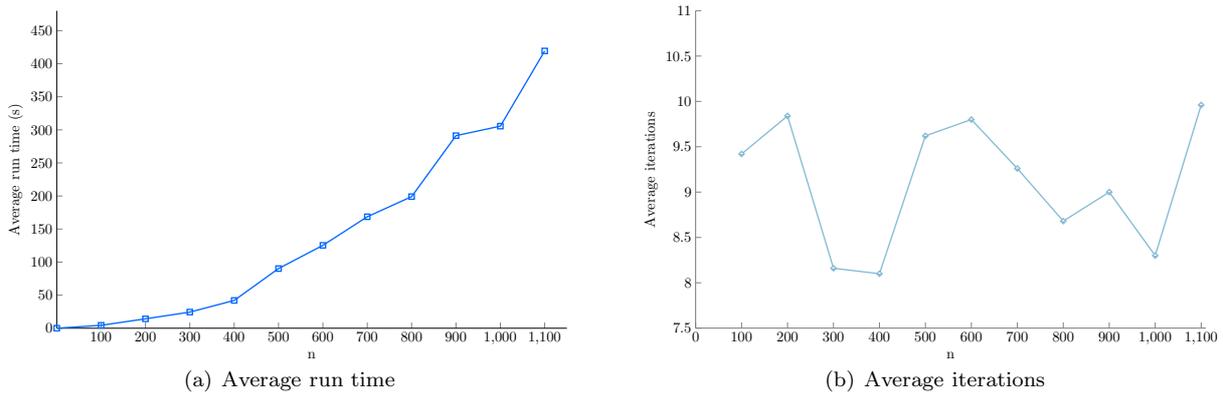

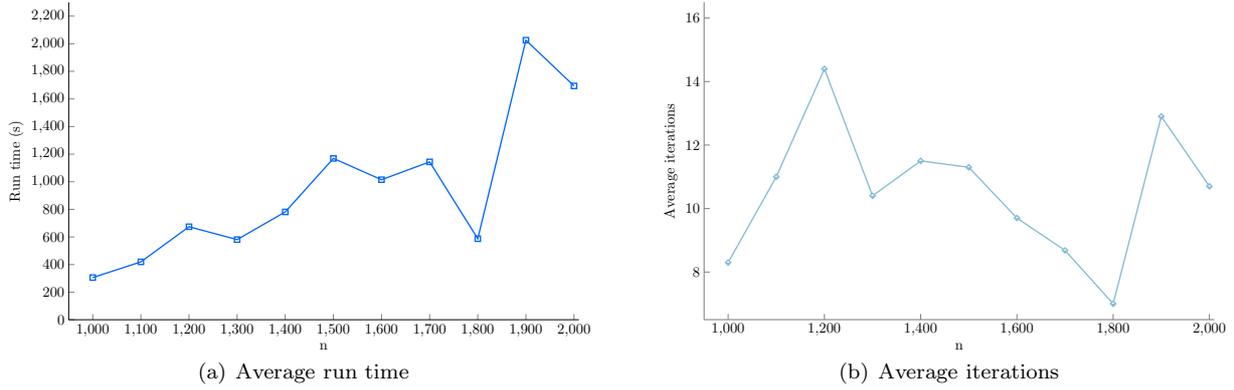
\begin{figure*}[h]
	\centering
	\subfigure[Average run time]{\makebox[6.5cm][c]{
			\begin{tikzpicture}[scale = 0.5]
				\begin{scope}[local bounding box=plot 1]
				\begin{axis}
					[name=plotter,
					xlabel={n}, ylabel={Run time (s)}, xmin=950, xmax=2010,
					ymin=0, ymax=2300,
					height=10cm,
					width=15cm,
					legend columns=5,
					legend pos=north east,
					axis y line*=left,
					axis x line*=bottom,
					ticklabel style = {font=\large},
					label style = {font=\large},
					xtick={1000,1100,1200,1300,1400,1500,1600,1700,1800,1900,2000},
					ytick={0,200,400, 600,800,1000,1200,1400,1600,1800,2000,2200},
					]
					
					\addplot
					[
					line width=1pt,
					color=blue!60!cyan,
					mark=square,
					]
					coordinates{
						(1000,305.57) (1100,419.30) (1200,673.9) (1300,580.69) (1400,781.19) (1500,1168.92) (1600,1014.67) (1700,1143.72) (1800,587.125) (1900,2025.71) (2000,1694.4)};
					
				\end{axis}
			\end{scope}
			\end{tikzpicture}
		}
	}
	\quad \quad \quad \quad\quad\quad
	\subfigure[Average iterations]{\makebox[6.5cm][c]{
			\begin{tikzpicture}[scale = 0.5]
			\begin{scope}[local bounding box=plot 2]
				\begin{axis}[
					name=zoom,
					width=50mm,
					height=50mm,
					at=(plotter.north east),
					axis line style={gray,thick},
					xlabel={n}, ylabel={Average iterations}, xmin=950, xmax=2010,
					ymin=6.5, ymax=16.5,
					height=10cm,
					width=15cm,
					legend columns=5,
					legend pos=north east,
					axis y line*=left,
					axis x line*=bottom,
					ticklabel style = {font=\large},
					label style = {font=\large},
					xtick={1000,1200,1400,1600,1800,2000},
					]
					\addplot
					[
					line width=1pt,
					color=cyan!50!gray!70,
					mark=halfsquare*,
					]
					coordinates{
						(1000,8.3) (1100,11) (1200,14.4) (1300,10.4) (1400,11.5) (1500,11.3) (1600,9.7) (1700,8.68) (1800,7.0) (1900,12.9) (2000,10.7)};
				\end{axis}
			\end{scope}
			\end{tikzpicture}
		}
	}
	\caption{Average run time and average iterations of the cutting plane method (for medium-size problems)}\label{fig-large}
\end{figure*}

The test instances were generated as follows:
\begin{enumerate}
	\item[1.] the size $n$ of the problem is between $100$ and $2000$;
	\item[2.] the integer $s$ is chosen randomly in $[1,10]$, and the vectors $v_1,\ldots,v_n\in \R^s$ are uniformly generated with each element $v_{ij}$ being in range $[1,10000]$ ($i=1,2,\ldots,n$, $j=1,\ldots,s$);
	\item[3.] $Q=[q_{ij}]$, where $q_{ij} = \norm{v_i - v_j}^2$ for all $i,j = 1,\ldots,n$;
	\item[4.] entries of $q=(q_1,\ldots,q_n)$ are randomly chosen within the range $[1,10000]$;
	\item[5.] the capacity $m$ is chosen randomly in $\{1,\ldots,n\}$.
\end{enumerate}
For medium size problems ($n=100,\ldots,1100$), we generated $50$ test instances for each dimension; and for large size problems ($n=1000,1100,\ldots, 2000$), we generated $10$ test instances for each dimension. Algorithm~\ref{Al-1} was used to solve each instance with the restriction on the number of iterations as $20$. 
The numerical experiments were
performed on a Dell Intel Core i7-8565U CPU 1.80GHz 16.0 GB, with the linear models solved
using the standard branch-and-cut optimizers in CPLEX version 12.10. The numerical results for medium size problems are presented in Table~\ref{100-1000} and Figure~\ref{fig-100-1000}. The numerical results for large size problems are presented in Table~\ref{larges} and Figure~\ref{fig-large}. The optimality gap in these tables is calculated as $\dfrac{\text{UB} - \text{LB}}{\text{UB}}\times100 \%$. The results show that the cutting plane method can derive high quality solutions for \eqref{QKP0} up to $n=2000$ in reasonable time.

\begin{figure*}[h]
	\centering
	\subfigure[Average run time]{\makebox[6.5cm][c]{
				\begin{tikzpicture}[scale = 0.5]
				\begin{scope}[local bounding box=plot 1]
					\begin{axis}
						[name=plotter,
						xlabel={n}, ylabel={Run time (s)}, xmin=50, xmax=100,
						ymin=0, ymax=240,
						height=10cm,
						width=15cm,
						legend columns=5,
						legend pos=north east,
						axis y line*=left,
						axis x line*=bottom,
						ticklabel style = {font=\large},
						label style = {font=\large},
						xtick={50,60,70,80,90,100},
						ytick={20,40,60,80,100,120,140,160,180,200,220,240},
						]
						
						\addplot
						[
						line width=1pt,
						color=blue!60!cyan,
						mark=square,
						]
						coordinates{
						(50,2.22) (60,3.25) (70,4.18) (80,10.24) (90,6.08) (100,5.4)};
						

						\addplot
						[
						line width=1pt,
						color=red!40!gray,
						mark=triangle,
						]
						coordinates{
							(50,200) (60,167.532) (70,200) (80,200) (90,200) (100,200)};

						\addplot
						[
						line width=1pt,
						color=green!30!gray,
						mark=halfsquare*,
						]
						coordinates{ (50,162.6) (60,141.738) (70,141.6)(80,180.5) (90,200) (100,183.743)};
						
						\legend{Cutting plane method, CPLEX,Glover's linearization}
					\end{axis}
				\end{scope}
			\end{tikzpicture}
		}
	}
	\quad \quad \quad \quad\quad\quad
	\subfigure[Average optimality gap]{\makebox[6.5cm][c]{
				\begin{tikzpicture}[scale = 0.5]
				\begin{scope}[local bounding box=plot 2]
					\begin{axis}
						[name=plotter2,
						at=(plotter.north east),
						xlabel={n}, ylabel={Optimality gap (\%)}, xmin=50, xmax=100,
						ymin=0, ymax=100,
						height=10cm,
						width=15cm,
						legend columns=5,
						legend pos=north east,
						axis y line*=left,
						axis x line*=bottom,
						ticklabel style = {font=\large},
						label style = {font=\large},
						xtick={50,60,70,80,90,100},
						ytick={10,20,30,40,50,60,70,80,90,100},
						]
						
						\addplot
						[
						line width=1pt,
						color=blue!60!cyan,
						mark=square,
						]
						coordinates{
							 (50,0.0) (60,0.0) (70,0) (80,0) (90,0) (100,0)};
						

						\addplot
						[
						line width=1pt,
						color=red!40!gray,
						mark=triangle,
						]
						coordinates{
							(50,35.691) (60,27.851) (70,36.685) (80,49.19) (90,43) (100,42)};

						\addplot
						[
						line width=1pt,
						color=green!30!gray,
						mark=halfsquare*,
						]
						coordinates{ (50,8.777) (60,11.056) (70,13.89) (80,29.63) (90,21) (100,15)};
						
						\legend{Cutting plane method, CPLEX,Glover's linearization}
					\end{axis}
				\end{scope}
			\end{tikzpicture}
		}
	}
	\caption{Performance comparison between: cutting plane method, CPLEX, and Glover's linearization.}	\label{comparision}
\end{figure*}
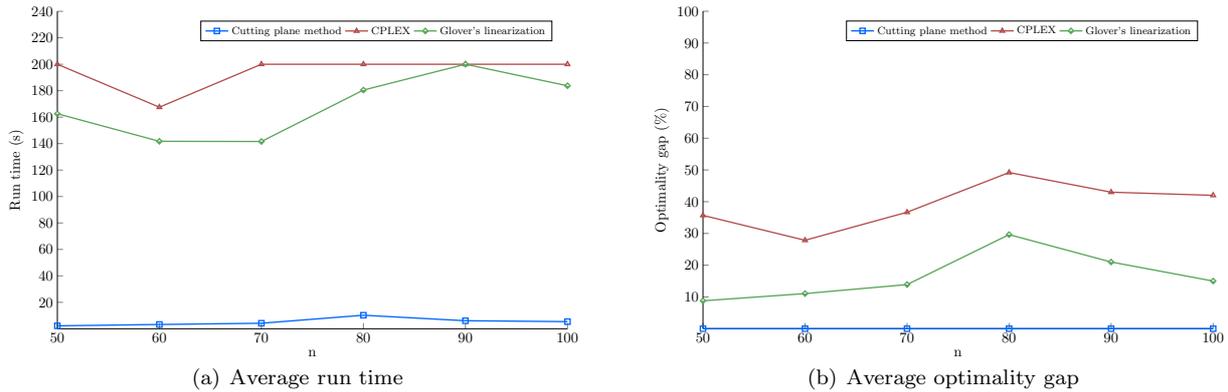

Algorithm~\ref{Al-1} was tested against Glover's linearization (see \cite{Glover,ADAMS200499}), and the mixed integer quadratic programming solver in CPLEX (see \cite{bliek1u2014solving}). Glover's linearization generally outperforms the standard linearization (see \cite{PISINGER2007623}), and therefore we only consider here the Glover's linearization method. However, both Glover's Linearization method and CPLEX cannot solve within $10\%$ optimality gap for instances with sizes more than $200$ in less than $200s$. Here, we only generate test instances for $n$ from $50$ to $100$ ($10$ instances for each dimension) and set the time limit as $200s$.  Table~\ref{TB3} and Figure~\ref{comparision} show that the cutting plane method still outperforms CPLEX and Glover's linearization for problems with $n$ between $50$ and $100$.

\begin{table}
\begin{center}
	\begin{tabular}{ccccc}
		\toprule
		$n$ & Average CPU time (s) & Average gap (\%) & Instances with zero gap& Average iterations\\ 
		\midrule
		$100$ & $4.44$ & $4.55\times 10^{-12}$ & $40/50$&  $9.42$ \\
		$200$ & $14.19$ & $1.26\times 10^{-3}$ & $42/50$& $9.84$\\  
		$300$ & $24.40$ & $3.59\times 10^{-4}$ &  $42/50$& $8.16$ \\
		$400$ & $42.00$ & $5.56\times 10^{-14}$& $45/50$&$8.10$ \\
		$500$ & $90.22$ & $1.59\times 10^{-4}$&$43/50$& $9.62$ \\
		$600$ & $125.34$ & $3.34\times 10^{-5}$ & $47/50$ &  $9.80$\\
		$700$ &  $168.55$ & $7.12\times 10^{-5}$ &  $45/50$ &$9.26$\\
		$800$ & $199.01$  &  $7.08\times 10^{-5}$  & $46/50$  &  $8.68$\\
		$900$ & $291.17$ & $5.27\times 10^{-5}$  &  $46/50$ & $10.30$ \\
		$1000$ & $305.57$ & $8.02\times 10^{-16}$  & $48/50$  & $9.00$  \\
		$1100$ & $419.30$  & $1.08\times 10^{-14}$ & $46/50$ & $9.96$\\
		\bottomrule
	\end{tabular}
\end{center}
\caption{Computational performance of the cutting plane method (for medium-size problems)}\label{100-1000}
\end{table}

\begin{table}
	\begin{center}
		\begin{tabular}{ ccccc}
				\toprule
			$n$ & Average CPU time (s) & Average gap (\%) & Instances with zero gap& Average iterations\\ 
					\midrule
			$1000$ & $226.69$ & $0.00$ & $10/10$&  $8.30$ \\
			$1100$ & $410.95$ & $0.00$ & $10/10$& $11.00$\\  
			$1200$ & $673.93$ & $3.93\times 10^{-4}$ &  $9/10$& $14.40$ \\
			$1300$ & $580.69$ & $3.39\times 10^{-15}$& $9/10$&$10.40$ \\
			$1400$ & $781.19$ & $1.83\times 10^{-14}$&$9/10$& $11.50$ \\
			$1500$ & $1168.92$ & $1.70\times 10^{-14}$ & $8/10$ &  $11.30$\\
			$1600$ &  $1014.66$ & $1.44\times 10^{-14}$ &  $9/10$ &$9.70$\\
			$1700$ & $1143.71$  &  $1.12\times 10^{-14}$  & $9/10$  &  $8.60$\\
			$1800$ & $587.12$ & $0.00$  &  $10/10$ & $7.00$ \\
			$1900$ & $2025.71$ & $1.10\times 10^{-5}$  & $8/10$  & $12.90$  \\
			$2000$ & $1694.43$  & $1.82\times 10^{-15}$ & $9/10$ & $10.70$\\
					\bottomrule
		\end{tabular}
	\end{center}
	\caption{Computational performance of cutting plane method (for large-size problems)}\label{larges}
\end{table}

\begin{table}
	\centering\begin{tabular}{cccccccc}
		\toprule
		&\multicolumn{3}{c}{Cutting plane method}&\multicolumn{2}{c}{IBM CPLEX  MIQP solver} & \multicolumn{2}{c}{Glover's linearization}\\
			\cmidrule(lr){2-4} \cmidrule(lr){5-6} \cmidrule(lr){7-8}
		$n$ & gap (\%) & run time (s) & iterations & gap (\%) & run time (s) &gap (\%) & run time (s)\\
		\midrule
		$50$ & $1.95\times10^{-12}$& $2.22$ & $9.60$ & $35.69$ & $200.00$ & $8.77$& $162.60$ \\
		$60$ & $1.43\times10^{-11}$& $3.25$ & $8.90$ & $27.85$ & $167.53$ & $11.05$ & $141.73$\\
		$70$ & $2.89\times10^{-15}$& $4.18$ & $8.90$ & $36.68$ & $200.00$ & $13.89$ & $141.62$\\
		$80$ & $1.69\times10^{-14}$& $10.24$ & $14.20$ & $49.19$ & $200.00$ & $29.63$ & $180.54$\\
		$90$ & $9.58\times10^{-13}$& $6.08$ & $9.50$ & $43.46$ & $200.00$ & $21.48$ & $200.00$\\
		$100$ & $4.21\times10^{-14}$& $5.40$ & $8.60$ & $42.23$ & $200.00$ & $16.06$ & $183.74$\\
		\bottomrule
	\end{tabular}
\caption{Comparison on average optimality gap and average run time between: cutting plane method, CPLEX, and Glover's linearization}\label{TB3}
\end{table}

\addcontentsline{toc}{section}{References}
\bibliographystyle{amsrefs}
\bibliography{Ref}
  \end{document}